\documentclass[12pt]{article}
 \usepackage{amssymb,amsmath,amsthm,amscd}
\usepackage{latexsym}
\usepackage{cite}
\usepackage{color} 

\newtheorem{theorem}{Theorem}[section]
\newtheorem{lemma}[theorem]{Lemma}

\newtheorem{definition}[theorem]{Definition}

\numberwithin{equation}{section}
\title{Entire solutions originating from three fronts to a two-dimensional nonlocal periodic lattice dynamical system}
\author{Shao-Hua Gan\footnote{ College of Science, University of Shanghai for Science and Technology,
Shanghai, 200093, China.
}\ \ and Zhixian Yu\footnote{The corresponding author. College of Science, University of Shanghai for Science and Technology,
Shanghai, 200093, China. Email: zxyu0902@163.com. Partially supported by Natural Science Foundation of Shanghai(No.18ZR1426500).
}
}
\date{}
\begin{document}
\maketitle
\begin{abstract}

This paper is concerned with the entire solutions of a two-dimensional nonlocal periodic lattice dynamical system.  With bistable assumption, it is well known that the system has three different types of traveling fronts. The existence of merging-front entire solutions originating from two fronts for the system
have been established by Dong, Li \& Zhang [{\it Comm. Pur Appl. Anal.},  {\bf17}(2018), 2517-2545]. Under certain conditions on the wave speeds, and by some auxiliary rational functions with certain properties to construct appropriate super- and sub solutions of the system, we establish two new types of entire solutions which originating from three fronts.

 \medskip

\textbf{AMS Subjective Classifications (2000): } 34K05; 34A34; 34E05
\end{abstract}

 \newpage

\section{Introduction}

\noindent

In this paper, we are interested in the entire solutions of the following  nonlocal periodic lattice dynamical system
\begin{equation}\label{eq1.1}
 u_{i,j}'(t)=\sum_{k_1\in\mathbb{Z}}\sum_{k_2\in\mathbb{Z}}J(k_1,k_2)u_{i-k_1,j-k_2}(t)-  u_{i,j}(t)+ f_{i,j}(u_{i,j}(t)),\,\,\, i, j \in \mathbb{Z},\,\, t\in \mathbb{R},
\end{equation}
 where $u_{i,j}(t)$ means the density of a certain species in a periodic patchy environment; $\sum\limits_{k_1\in\mathbb{Z}}\sum\limits_{k_2\in\mathbb{Z}}J(k_1,k_2)u_{i-k_1,j-k_2}(t)-u_{i,j}(t)$ is the nonlocal
 dispersal and represents transportation due to long range dispersion mechanism. The kernel function
 $J:\mathbb{Z}\times\mathbb{Z}\longrightarrow[0,\infty)$ is a probability function formulating the dispersal of individuals and
 satisfies
 $J(\cdot,\cdot)\geq0$ is even $\sum\limits_{k_1\in\mathbb{Z}}\sum\limits_{k_2\in\mathbb{Z}}J(k_1,k_2)=1$, $J(k_1,k_2)=0$ if $|k_1|>k_0$ or $|k_2|>k_0$, where $k_0$ is a
 positive constant.
Here the reaction term $f_{i,j}\in C^2 (\mathbb{R})$ satisfies the bistable condition
\begin{enumerate}
\item[$\rm(A1)$] $f_{i,j}(0)=f_{i,j}(a)=f_{i,j}(1)=0$, $i, j\in\mathbb{Z}$, $f_{i,j}'(0)<0$, $f_{i,j}'(1)<0$, $f_{i,j}'(a)>0$, $f_{i,j}(u)<0$\: for \:$u\in(0,a)$\: and \: $f_{i,j}(u)>0$\: for \: $u\in(a,1)$
\end{enumerate}
and the periodic condition
\begin{enumerate}
\item[$\rm(A2)$] $ f_{i,j}(\cdot)= f_{i+N_1,j}(\cdot)= f_{i,j+N_2}(\cdot)$ for all $N_1$ and $N_2$ are two positive integers.
\end{enumerate}

The system \eqref{eq1.1} can be also regarded as a nonlocal version of the following local diffusion system
\begin{equation}\label{eq1.2}
 u_{i,j}'(t)=D[u_{i,j}](t)- d_{i,j} u_{i,j}(t)+ f_{i,j}(u_{i,j}(t)),
\end{equation}
where
\begin{align*}
D[u_{i,j}](t)&=p_{i+1,j} u_{i+1,j}(t)+p_{i,j}u_{i-1,j}(t)+q_{i,j+1}u_{i,j+1}(t)+q_{i,j}u_{i,j-1}(t),\\
d_{i,j}&=p_{i+1,j}+p_{i,j}+q_{i,j+1}+q_{i,j}.
\end{align*}
where $ p_{i,j},q_{i,j}> 0$ and $p_{i+N_1,j} = p_{i,j} = p_{i,j+N_2}$, $q_{i+N_1,j} = q_{i,j} =
q_{i,j+N_2}$. The reaction function satisfies $f_{i,j}(\cdot) = f_{i+N_1,j}(\cdot) = f_{i,j+N_2}(\cdot)$ for all $ i,j \in {\mathbb Z}$.

In the past decades, there are a lot of works devoted to front propagation for lattice dynamical systems in spatially periodic or homogeneous media, e.g.\cite{PWAC,CPWT,CWG,CCWU, CPYH,CPWTZC,CPCWT,CPCWTZ,JACU,XCSC,XCJSU,JSGFH,JSGCHW,JSGCHWE,JSGCHWT,JSGYW,SMPW} and some references cited therein. Notice that the traveling front solution is a special type of entire solutions, but it is not enough to
understand the whole dynamics. From the viewpoint of dynamics,  for entire solutions which behave as
two traveling wave fronts moving towards each other from both sides of the $x$ (or $j$) axis, this type of entire solutions are called ``{\it
annihilating-front}'' entire solutions. Beyond that, there are also two common types of entire solutions. The first type behaves as a
monostable front merges with a bistable front and one chases another from the same side; while the other type can be represented by
two monostable fronts approaching each other from both sides of the $x$ (or $j$) axis and merging and converging to a single bistable front, see \cite{wushiyang,Wang7,CHWU,ZhangLiWu,Moritan}. Such two types of entire solutions are often called as merging-front entire solutions. In \cite{Dongli}, the authors have established the existence of merging-front entire solutions originating from two fronts for the system of (\ref{eq1.1}) with both monostable and bistable nonlinearities.\medskip

Recently, there are other new types of entire solutions merging three fronts, which were addressed in \cite{ChenGuoNinomiyaYao,YYC}. Motivated by these works, it is natural and interesting to study new entire solutions merging three fronts of the system (\ref{eq1.1}).
Now, there has been no results on the entire solutions merging three fronts for the nonlocal periodic lattice dynamical system.
Therefore, in this paper, we will study the existence of the entire solutions originating from three fronts of \eqref{eq1.1}.\medskip

Before to state our main results, we first give definitions of the pulsating traveling fronts and entire solutions for \eqref{eq1.1}.
\begin{definition}[Pulsating traveling fronts and entire solutions]\label{def1.1}~
\begin{enumerate}
\item[{\rm (1)}] A solution $u(t)=\{u_{i,j}(t)\}_{i,j\in\mathbb{Z}}$, $t\in\mathbb{R}$ of \eqref{eq1.1} is called a pulsating
{\rm(}or periodic{\rm)} traveling front connecting $\{e_1,e_2\}\subset\{0,a,1\}$ in the direction $(cos\theta$, $sin\theta)$ with the wave speed $k$,
if
\begin{eqnarray}
u_{i,j}(t)=\phi_{i,j}(icos\theta +jsin\theta+kt),  \text{ or }
u_{i,j}(t)=\phi_{i,j}(-icos\theta-jsin\theta+kt),\nonumber
\end{eqnarray}
for all $i, j\in\mathbb{Z},{\ }t\in \mathbb{R}$ and some function $\phi(\cdot)=\{\phi_{i,j}(\cdot)\}_{i,j\in\mathbb{Z}}$ which satisfies
\begin{eqnarray}
\phi_{i,j}(\cdot)= \phi_{i+N_1,j}(\cdot)= \phi_{i,j+N_2}(\cdot),\ \phi_{i,j}(-\infty)=e_1 \text{ and  }\phi_{i,j}(\infty)=e_2,{\ \ }\forall i,j\in\mathbb{Z}.\nonumber
\end{eqnarray}
\item[{\rm (2)}] A function $u(t)=\{u_{i,j}(t)\}_{i,j\in {\mathbb Z}}$, $t\in\mathbb{R}$ is called an entire solution of \eqref{eq1.1} if
for any $i,j\in{\mathbb Z}$, $u_{i,j}(t)$  is differentiable for all $t\in{\mathbb R}$ and satisfies \eqref{eq1.1} for $i,j\in{\mathbb Z}$ and $t\in{\mathbb R}$.
\end{enumerate}
\end{definition}

According to  \cite[Theorem 6]{ChenGuo}, it is easy to see that \eqref{eq1.1} has two increasing pulsating traveling fronts
$\tilde\varphi_{i,j;1}(icos\theta +jsin\theta+\tilde{v}t)$ and $\varphi_{i,j;1}(-i cos\theta-j sin\theta+\bar{v}t)$ satisfying
\begin{eqnarray}
\tilde\varphi_{i,j;1}(-\infty)=\varphi_{i,j;1}(-\infty)=0 \text{ and  } \tilde\varphi_{i,j;1}(+\infty)=\varphi_{i,j;1}(+\infty)=1,{\ \ } \forall i,j\in\mathbb{Z},
\label{1.5}
\end{eqnarray}
where $\tilde{v}$ and $\bar{v}>0$.
Let $\hat{v}=-\bar{v}$, then $\hat\varphi_{i,j;1}(icos\theta +jsin\theta+\hat{v}t)=\varphi_{i,j;1}(-i cos\theta-j sin\theta+\bar{v}t)=\varphi_{i,j;1}(-(icos\theta +jsin\theta+\hat{v}t))$, which satisfies
\begin{eqnarray}
\hat\varphi_{i,j;1}(-\infty)=1 \text{ and  }\hat\varphi_{i,j;1}(+\infty)=0,{\ \ } \forall i,j\in\mathbb{Z}.\label{1.6}
\end{eqnarray}

Moreover, if we restrict $f_{i,j}(u)$ in the intervals  $[0,a]$ and $[a,1]$, respectively, then \eqref{eq1.1} can be regarded as two
monostable equations. Thus, we make the following assumption:

\noindent \ \ (A3)\quad
$f_{i,j}(u)\Big\{
\begin{array}{ll}
\geq f'_{i,j}(a)(u-a),&\mbox{if }u\in[0,a],\\
\leq f'_{i,j}(a)(u-a),&\mbox{if }u\in[a,1],
 \end{array}\mbox{for all } i,j\in\mathbb{Z}.$\medskip

For the case $u_{i,j}(t)\in[0,a]$, by considering $ v_{i,j}(t)=a-u_{i,j}(t)$, \eqref{eq1.1} can be reduce to
\begin{equation}\label{eq1.5}
 v_{i,j}'(t)=\sum_{k_1\in \mathbb{Z}}\sum_{k_2\in\mathbb{Z}}J(k_1,k_2)v_{i-k_1,j-k_2}(t)-v_{i,j}(t)+ g_{i,j}(v_{i,j}(t))\,\,\mbox{for any}\,\, i, j \in \mathbb{Z},\,\, t\in \mathbb{R},
\end{equation}
where $g_{i,j}(v)=-f_{i,j}(a-v)$. Then it follows from \cite{ZhangLiWu,Dongli} that there exists a $v^{*}(\theta)>0$ such that for every $\check{v}_1, \hat{v}_1>v^{*}$, \eqref{eq1.5} has two increasing pulsating traveling fronts $\check{W}_{i,j}(icos\theta +jsin\theta+\check{v}_1t)$ and $\hat{W}_{i,j}(icos\theta +jsin\theta+\hat{v}_1t)$ which satisfy
\begin{center}
$
  \check{W}_{i,j}(-\infty)=\hat{W}_{i,j}(-\infty)=0 \text{ and  }\check{W}_{i,j}(+\infty)=\hat{W}_{i,j}(+\infty)=a,
$  for $i,j\in\mathbb{Z}$.
\end{center}
Let's define $v^*_1=-v^*$, $\tilde{v}_1=-\check{v}_1<-v^*=v^*_1$,
$\tilde{\varphi}_{i,j;2}(icos\theta+jsin\theta+\tilde{v}_1t)=a-\check{W}_{i,j}(-icos\theta-jsin\theta+\check{v}_1t)
=a-\check{W}_{i,j}(-(icos\theta+jsin\theta)+\tilde{v}_1t)$ for any $\check{v}_1<v^*_1$\\
and
$\hat{\varphi}_{i,j;2}(icos\theta+jsin\theta+\hat{v}_1t)=a-\check{W}_{i,j}(icos\theta+jsin\theta+\hat{v}t)$ for any $\hat{v}_1>v^{*}$ then $\tilde{\varphi}_{i,j;2}(icos\theta+jsin\theta+\tilde{v}_1t)$ and $\hat{\varphi}_{i,j;2}(icos\theta+jsin\theta+\hat{v}_1t)$ are pulsating traveling fronts of \eqref{eq1.1} with the wave speeds $\tilde{v}_1<v^*_1$ and $\hat{v}_1>v^{*}$, respectively,
\begin{center}
$\tilde{\varphi}_{i,j;2}(-\infty)=0$, $\tilde{\varphi}_{i,j;2}(+\infty)=a$, and $\hat{\varphi}_{i,j;2} (-\infty)=a$, $\hat{\varphi}_{i,j;2}(+\infty)=0$
\end{center}
Similarly, for the case $u_{i,j}(t)\in[a,1]$, we can also obtain that (letting $v_{i,j}(t)=u_{i,j}(t)-a)$ there exists a ${v_2}^{*}(\theta)>0$
such that for every $v_2>{v_2}^{*}$, \eqref{eq1.1} admits an increasing pulsating traveling front $\tilde\varphi_{i,j;3}(j+v_2t)$
with
\begin{center}
$\tilde\varphi_{i,j;3}(-\infty)=a$ and $\tilde\varphi_{i,j;3}(+\infty)=1$ for any $icos\theta+jsin\theta\in\mathbb{Z}$.
\end{center}

For \eqref{eq1.1}, the traveling front $u_{i,j}(t)$  with the speed $v$ exists if $u_{i,j}(t)=\varphi(icos\theta+jsin\theta+vt)$ for
some function and it connects two different two different constant states. Now we set $z:=i cos\theta +j sin\theta+vt$ and substitute
into \eqref{eq1.1}, then
let $(\varphi_{i,j},v_i),i=1,2,3$ be the traveling fronts of \eqref{eq1.1}. If the entire solution $ \varphi_{i,j}(t)$ of
\eqref{eq1.1} satisfies
\begin{align*}
&\ \lim\limits_{t\rightarrow-\infty}\Big\{ \sum_{1\leq i\leq3}\sup_{ p_{i-1}(t)<i cos\theta +j sin\theta<p_i(t)}
|W_{i,j}(t)-\varphi_{i,j}(t)(i cos\theta +j sin\theta+v_{i}t+\theta_{i})|\Big\}=0,
\end{align*}
where $v_1<v_2<v_3,$ $\theta_{1},\theta_{2},\theta_{3}$ are some constants, $p_i(t):=-(v_i+v_{i+1})t/2$, $p_0=-\infty$ and $p_3=+\infty$, it is called the
entire solution originating from three fronts of \eqref{eq1.1}.

\begin{theorem}\label{thm1.2}
Assume that (A1), (A2) and (A3) hold. Let $(\hat{v},\hat\varphi_{i,j;1}),\,\, (\tilde v_1,\tilde\varphi_{i,j;2})$ and $ (v_2,\tilde\varphi_{i,j;3})$ be pulsating traveling fronts described as above such that
$\hat{v}<\tilde v_1<{v_2}$. Then there exists an entire solution $U_{i,j}(t):\mathbb{Z}\times\mathbb{R} \rightarrow
[0,1]$ originating from three fronts of \eqref{eq1.1} and $\omega\in \mathbb{R}$ which satisfy
$$\lim\limits_{t\rightarrow+\infty}\sup_{i cos\theta +j sin\theta\in\mathbb{R}}
|U_{i,j}(t)-1|=0$$
and
\begin{align}
&\lim_{t\rightarrow-\infty}\big\{\sup_{i cos\theta +j sin\theta\leq m_1 }
|U_{i,j}(t)-\hat{\varphi}_{i,j;1}(i cos\theta +j sin\theta+\hat{v}t-\omega)|\nonumber\\
&\qquad\qquad\qquad\nonumber\\
&\qquad\qquad\qquad+\sup_{m_1 \leq i cos\theta +j sin\theta\leq m_2}
|U_{i,j}(t)-\tilde{\varphi}_{i,j;2}(i cos\theta +j sin\theta+\tilde v_1t+\omega)|\nonumber\\
&\qquad\qquad\qquad\nonumber\\
&\qquad\qquad\qquad+\sup_{i cos\theta +j sin\theta\geq m_2}
|U_{i,j}(t)-\tilde{\varphi}_{i,j;3}(i cos\theta +j sin\theta+v_2t+\omega)|\big\}=0,\nonumber
\end{align}
where
$m_1=-\frac{(\hat {v}+\tilde v_1)t}{2}$ and $m_2=-\frac{(\tilde v_1+v_2)t}{2}$.

\end{theorem}

\begin{theorem}\label{thm1.3}  Assume that (A1), (A2) and (A3) hold. Let $(\hat{v},\hat\varphi_{i,j;1}),$ $(\tilde v_1,\tilde\varphi_{i,j;2})$ and $(\hat v_1,\hat\varphi_{i,j;2})$ be pulsating traveling fronts with
$\hat{v}<\tilde v_1<\hat v_1$. Then there exists an entire solution $V_{i,j}(t):\mathbb{Z}\times\mathbb{R}\rightarrow[0,1]$ originating from three fronts of \eqref{eq1.1} and $\omega_1,\omega_2\in\mathbb{R}$
which satisfy
$$\lim\limits_{t\rightarrow+\infty}\sup_{i cos\theta +j sin\theta\in\mathbb{R}}|V_{i,j}(t)- \hat{\varphi}_{i,j;1}(i cos\theta +j sin\theta+\hat{v}t+\omega_1)|=0$$
and
\begin{align*}
&\lim_{t\rightarrow-\infty}\big\{\sup_{i cos\theta +j sin\theta\leq n_1}
|V_{i,j}(t)-\hat{\varphi}_{i,j;1}(i cos\theta +j sin\theta+\hat{v}t+\omega_1)|\nonumber\\
&\qquad\qquad\qquad\nonumber\\
&\qquad\qquad\qquad+\sup_{n_1 \leq i cos\theta +j sin\theta\leq n_2}
|V_{i,j}(t)-\tilde{\varphi}_{i,j;2}(i cos\theta +j sin\theta+\tilde v_1t+\omega_1)|\nonumber\\
&\qquad\qquad\qquad\nonumber\\
&\qquad\qquad\qquad+\sup_{i cos\theta +j sin\theta\geq n_2}
|V_{i,j}(t)-\hat{\varphi}_{i,j;2}(i cos\theta +j sin\theta+\hat v_1t-\omega_2)|\big\}=0,
\end{align*}
where $n_1=-\frac{(\hat {v}+\tilde v_1)t}{2}$ and $n_2=-\frac{(\tilde v_1+\hat v_1)t}{2}$.

\end{theorem}

In order to verify Theorems \ref{thm1.2}-\ref{thm1.3}, we will adopt the elementary method of super- and subsolutions and a comparison principle. Since the comparison principle can be well applied, we only need to construct a suitable pair of super- and subsolutions by some auxiliary rational functions with certain properties which were developed by Morita and Ninomiya in \cite{Moritan}. This technique had been used to prove some types of entire solutions originating from two fronts of \eqref{eq1.1}. Therefore, we would apply the technique to establish some new types of entire solutions originating from three fronts of \eqref{eq1.1}, i.e. Theorems \ref{thm1.2}-\ref{thm1.3}.

Now we recall the definitions of super- and subsolutions of \eqref{eq1.1} and some known results on the existence, the comparison principle and the prior estimates for solutions of \eqref{eq1.1}.
\begin{definition}[Super- and subsolutions]\label{def2.1} ~
\begin{enumerate}
\item[{\rm(1)}] Let $t_0<T$ be any two real constants.
A sequence of continuous functions $\{u_{i,j}(t)\}_{i,j\in {\mathbb Z}}$
is called a supersolution {\rm(}or subsolution{\rm)} of \eqref{eq1.1} on $[t_0,T)$ if
\begin{eqnarray*}
u_{i,j}'(t)\geq (or \leq)\ \sum_{k_1\in\mathbb{Z}}\sum_{k_2\in\mathbb{Z}}J(k_1,k_2)u_{i-k_1,j-k_2}(t)-  u_{i,j}(t)+ f_{i,j}(u_{i,j}(t)),
\end{eqnarray*}
for $t\in [t_0,T)$.
\item[{\rm(2)}] Let $T$ be a real constant. A sequence of continuous functions $\{u_{i,j}(t)\}_{i,j\in {\mathbb Z}}$
is called a supersolution {\rm(}or subsolution{\rm )} of \eqref{eq1.1} on $(-\infty,T)$ if for any $t_0<T$, $\{u_{i,j}(t)\}_{i,j\in {\mathbb Z}}$ is a
supersolution {\rm (}or subsolution{\rm)} of \eqref{eq1.1} on $[t_0,T)$.
\end{enumerate}
\end{definition}

\begin{lemma}\cite[Existence and Comparison principle]{Dongli}\label{com}~
\begin{enumerate}
\item[{\rm(1)}] For any $\varphi=\{\varphi_{i,j}\}_{i,j\in {\mathbb Z}}$ with $\varphi_{i,j}\in[0,1]$,
equation \eqref{eq1.1} admits a unique solution $u(t;\varphi)=\{u_{j}(t;\varphi)\}_{i,j\in {\mathbb Z}}   $ on $[t_0,+\infty)$ which
satisfies $u_{i,j}(t_0;\varphi)=\varphi_{i,j}$ and $u_{i,j}\in C^1([t_0,\infty),[0,1])$ for all $i,j\in {\mathbb Z}$.
\item[{\rm(2)}]
Suppose $\big\{u^+_{i,j}(t)\big\}_{i,j\in {\mathbb Z}}$ and $\big\{u^-_{i,j}(t)\big\}_{i,j\in {\mathbb Z}}$
are bounded supersolution and subsolution of \eqref{eq1.1} on $[t_0,+\infty)$, respectively, such that $ u^+_{i,j}(t), u^-_{i,j}(t)\geq0$ and $u^+_{i,j}(t_0)\geq
u^-_{i,j}(t_0)$ for $i,j\in {\mathbb Z}$ and $t\geq t_0$, then $u^+_{i,j}(t)\geq
u^-_{i,j}(t)\geq0$ for all $i,j\in {\mathbb Z}$ and $t\geq t_0$.
\end{enumerate}
\end{lemma}
\begin{lemma}\cite[A Prior estimate]{Dongli}\label{com}~
Let $u(t;\varphi)=\big\{u_{i,j}(t;\varphi)\big\}_{i,j\in \mathbb{Z}}$ be a solution of \eqref{eq1.1} with initial value
$\varphi=\{\varphi_{i,j}\}_{i,j\in \mathbb{Z}}$ and $\varphi_{i,j}\in [0,1]$. For any fixed $t_0\in\mathbb{R}$, there exists a
positive constant $M$ {\rm(}independent of $ t_0$ and $\varphi${\rm)} such that
\begin{align*}
\big| u_{i,j}'(t;\varphi)\big|\leq M \text{ and } \big| u_{i,j}''(t;\varphi)\big|\leq M, \ \mbox{for }
{i,j}\in \mathbb{Z}\ \mbox{and }t> t_0+1.
\end{align*}
\end{lemma}

\section{Proof of Theorem 1.2}

Set $c_1:=\hat{v}$, $c_2:=\tilde v_1$ and $c_3:=v_2$. let
$\phi_{i,j;1}=\hat{\varphi}_{i,j;1}$, $\phi_{i,j;2}=\tilde{\varphi}_{i,j;2}$, $\phi_{i,j;3}=\tilde{\varphi}_{i,j;3}$, $i,j\in \mathbb{Z}, k=1,2,3,$
 be traveling fronts of  \eqref{eq1.1} that satisfy

\begin{equation}\label{wave}
\left\{
\begin{array}{ll}
{c_k}\phi_{i,j;k}'(\xi)=\sum_{k_1}\sum_{k_2}J(k_1,k_2)\phi_{i-k_1,j-k_2;k}(\xi)-\phi_{i,j;k}(\xi)+f_{i,j}(\phi_{i,j;k}(\xi)),\medskip\\
\phi_{i+N_1,j;k}(\xi)=\phi_{i,j;k}(\xi)=\phi_{i,j+N_2;k}(\xi),    i,j\in \mathbb{Z},~\xi\in \mathbb{R},\\
 \phi_{i,j;k}(-\infty)=\alpha_i, \phi_{i,j;k}(\infty)=\beta_i
\end{array}\right.
 \end{equation}
where $(\alpha_1,\beta_1,\alpha_2,\beta_2,\alpha_3,\beta_3)=(1,0,0,a,a,1)$.

Without loss of generality, we assume that
\begin{align}\label{2.2}
 \phi_{i,j;1}(0)=\frac{a}{2}, \phi_{i,j;2}(0)=\frac{a}{2}, \phi_{i,j;3}(0)=\frac{1+a}{2}.
\end{align}


\begin{lemma}\label{lem2.1} \
According to the arguments in \rm{\cite[Lemma 2.5]{Dongli}}, there exist positive numbers $C_{0},\, C_{1},$ $ C_{2}$, $\eta_{1},\, \eta_{2}$ and $\rho$, which depend on $\hat{v}$, $\tilde v_1$, $v_2$, $\hat{v}$, $\tilde v_1$, $\hat v_1$ and $\theta$, such that for all $ i,j\in\mathbb{Z}$ and $\xi\leq 1$,
\begin{align}
&0\leq \phi_{i,j;1}(t)'(\xi),\phi_{i,j;2}(t)'(\xi),\phi_{i,j;3}(t)'(\xi)  \leq C_{0}e^{\eta_{1}\xi},\label{a1}\medskip\\
C_{1}e^{\eta_{1}\xi}&\leq 1-\phi_{i,j;1}(t)(\xi),\phi_{i,j;2}(t)(\xi),\phi_{i,j;3}(t)(\xi)-a \leq
C_{2}e^{\eta_{1}\xi},\label{a3}\medskip\\
&\frac{\phi_{i,j;1}(t)'(\xi)}{ 1-\phi_{i,j;1}(t)(\xi)},\ \frac{\phi_{i,j;2}(t)'(\xi)}{\phi_{i,j;2}(t)(\xi)},\
\frac{\phi_{i,j;3}(t)'(\xi)}{\phi_{i,j;3}(t)(\xi)-a}\geq\rho; \label{a5}
\end{align}
For all $ i,j\in\mathbb{Z}$ and $\xi\geq -1$,
\begin{align}
&0\leq \phi_{i,j;1}(t)'(\xi),\phi_{i,j;2}(t)'(\xi),\phi_{i,j;3}(t)'(\xi)\leq C_{0}e^{-\eta_{2}\xi},\label{a2}\medskip\\
C_{1}e^{-\eta_{2}\xi}&\leq \phi_{i,j;1}(t)(\xi),a-\phi_{i,j;2}(t)(\xi),1-\phi_{i,j;3}(t)(\xi)\leq C_{2}e^{-\eta_{2}\xi},
\label{a4}\medskip\\
&\frac{\phi_{i,j;1}(t)'(\xi)}{ \phi_{i,j;1}(t)(\xi)},\ \frac{\phi_{i,j;2}(t)'(\xi)}{a-\phi_{i,j;2}(t)(\xi)},\
\frac{\phi_{i,j;3}(t)'(\xi)}{1-\phi_{i,j;3}(t)(\xi)}\geq\rho.\label{a6}
\end{align}
\end{lemma}

\subsection{Construction of sub- and supersolutions}
To construct a pair of super-solution and sub-solution of \eqref{eq1.1}, we introduce the following auxiliary function (see \cite{ChenGuoNinomiyaYao})
\begin{align}
Q(y,z,w):=&z+(1-z)\frac{(1-y)z(w-a)+y(a-z)(1-w)}{(1-y)z(1-a)+(a-z)(1-w)}\ \ \mbox{for any}\,\,(y,z,w)\in D_{1},\label{2.9}
\end{align}
where $D_{1}:=\{[0,1]\times[0,a]\times[a,1]\}\backslash(\{(1,a,w)| a \leq w\leq1\}\cup\{(1,z,1)| 0 \leq z \leq a \}\cup\{(y,0,1)|0\leq y\leq 1\}) $.
\begin{lemma}\cite[Lemma 3.1]{ChenGuoNinomiyaYao}\label{lem2.2} The function $Q(y,z,w)$ has the following properties.
\begin{itemize}
  \item [(i)]
\begin{equation}
Q(y,z,w)=\left\{ {\begin{array}{ll}
\displaystyle y+(1-y)z\frac{(1-a)(w-y)}{(1-y)z(1-a)
+(a-z)(1-w)}, \medskip\\
\displaystyle w+( a-z )(1-w)\frac{y-w}{(1-y)z(1-a)
+(a-z)(1-w)}. \\
\end{array}} \right.
\label{2.10}
\end{equation}
  \item [(ii)] There exist functions $Q_i,i=1,2,3$, such that
\begin{align*}
&Q_y(y,z,w)=(a-z)(1-w)Q_1(y,z,w),\\
&Q_z(y,z,w)=(1-y)(1-w)Q_2(y,z,w),\\
&Q_w(y,z,w)=(1-y)z Q_3(y,z,w).\\
\end{align*}
  \item [(iii)] There exist functions $R_j,j=1,\cdot\cdot\cdot,16,$ such that
\begin{align*}
Q_{yy}(y,z,w)=&z R_1(y,z,w)=(a-z) R_2(y,z,w)=(1-w)R_3(y,z,w),\\
Q_{zz}(y,z,w)=&(1-y) R_4(y,z,w)=(1-w) R_5(y,z,w)\\
=&y R_6(y,z,w)+(w-a)R_7(y,z,w),\\
Q_{ww}(y,z,w)=&(1-y) R_8(y,z,w)=z R_9(y,z,w)=(a-z)R_{10}(y,z,w),\\
Q_{yz}(y,z,w)=&(1-w) R_{11}(y,z,w),Q_{zw}(y,z,w)=(1-y) R_{12}(y,z,w),\\
Q_{yw}(y,z,w)=&(1-y) R_{13}(y,z,w)=z R_{14}(y,z,w)\\
=&(a-z)R_{15}(y,z,w)+(1-w)R_{16}(y,z,w).\\
\end{align*}
\end{itemize}
\end{lemma}

Next, we consider the following four ordinary differential equations with the initial conditions (see \cite{ChenGuoNinomiyaYao,Guom,Moritan}):
\begin{align}
&p_{1}'(t)=s_1+Le^{\kappa p_{1}(t)},\,\,p_{1}(0)=p_0;
\label{eq2.11}\\
&p_{2}'(t)=s_2+Le^{\kappa p_{1}(t)},\,\,p_{2}(0)=p_0;
\label{eq2.12}\\
&r_{1}'(t)=s_1-Le^{\kappa r_{1}(t)},\,\,r_{1}(0)=r_0;
\label{eq2.13}\\
&r_{2}'(t)=s_2-Le^{\kappa r_{1}(t)}, \,\,r_{2}(0)=r_0,
\label{eq2.14}
\end{align}
where $L,s_2,s_1,p_0,r_0$ and $\kappa$ are positive constants.
\eqref{eq2.11}-\eqref{eq2.14} have the following unique solutions for $t\leq 0$, respectively,
\begin{equation}
\label{eq2.15}
p_{1}(t)=s_1t-\frac{1}{\kappa}\ln\big[e^{-\kappa p_{0}}+\frac{L(1-e^{\kappa s_1 t})}{s_1}\big],
\end{equation}
\begin{equation}
\label{eq2.16}
p_{2}(t)=s_2t-\frac{1}{\kappa}\ln\big[e^{-\kappa p_{0}}+\frac{L(1-e^{\kappa s_1 t})}{s_1}\big],
\end{equation}
\begin{equation}
\label{eq2.17}
r_{1}(t)=s_1t-\frac{1}{\kappa}\ln\big[e^{-\kappa r_{0}}+\frac{L(1-e^{\kappa s_1 t})}{s_1}\big],
\end{equation}
\begin{equation}
\label{eq2.18}
r_{2}(t)=s_2t-\frac{1}{\kappa}\ln\big[e^{-\kappa r_{0}}+\frac{L(1-e^{\kappa s_1 t})}{s_1}\big].
\end{equation}
For a sufficiently large constant $\delta>0$, we choose $p_0$ and $r_0$ satisfying
\begin{align}\label{eq2.19}
p_0=-\frac{1}{\kappa}\ln(e^{-\kappa r_{0}}-\frac{2L}{s_1})<-\delta,\,\,r_0<-\frac{1}{\kappa}\ln(\frac{2L}{s_1}+e^{\kappa\delta}).
\end{align}
Moreover,
\begin{align}
&\lim\limits_{t\rightarrow-\infty}(p_{1}(t)-r_{1}(t))=\lim\limits_{t\rightarrow-\infty}(p_{2}(t)-r_{2}(t))=0,
\label{eq2.20}\\
&\lim\limits_{t\rightarrow-\infty}(p_{1}(t)-s_{1}t)=\lim\limits_{t\rightarrow-\infty}(p_{2}(t)-s_{2}t)=-\frac{1}{\kappa}\ln\big[e^{-\kappa p_{0}}+\frac{L}{s_1}\big],
\label{eq2.21}\\
&\lim\limits_{t\rightarrow-\infty}(r_{1}(t)-s_{1}t)=\lim\limits_{t\rightarrow-\infty}(r_{2}(t)-s_{2}t)=-\frac{1}{\kappa}\ln\big[e^{-\kappa r_{0}}+\frac{L}{s_1}\big].
\label{eq2.22}
\end{align}
Also, there exist a positive constant $N$ such that for all $t\leq 0$
\begin{equation}
\label{eq2.23}
0<p_1 (t)-r_1 (t)=p_2 (t)-r_2 (t)\le Ne^{ks_1 t},
\end{equation}
and $p_{1}(t), p_{2}(t), r_{1}(t), r_{2}(t)\leq-\delta$.

Therefore, by the choice of $k$, there exists a $t_0 <0$ such that
\begin{equation}
\label{eq2.24}
\max \big\{\frac{\eta _1 (p_2 (t)-p_1 (t))}{2},\frac{\eta _2 (p_2 (t)-p_1
(t))}{2}\big\}<kp_1 (t)<0,
\end{equation}
\begin{equation}
\label{eq2.25}
\max \big\{\frac{\eta _1 (r_2 (t)-r_1 (t))}{2},\frac{\eta _2 (r_2 (t)-r_1
(t))}{2}\big\}<kr_1 (t)<0,\ \forall t\le t_0.
\end{equation}

\begin{lemma}\label{lem2.3} \
Let $p_{1}(t)$ and $ p_{2}(t)$ be solution of \eqref{eq2.15} and \eqref{eq2.16}. For any $p_{2}(t)\leq p_{1}(t)\leq 0$. There exist positive constants $\epsilon_1$,
$\epsilon_2$ and $\epsilon_3$ such that
\begin{align*}
&Q_y(\phi_{i,j;1}(z-p_{1}(t)),\phi_{i,j;2}(z+p_{1}(t)),\phi_{i,j;3}(z+p_{2}(t)))\geq \varepsilon_1 ,\ \ \mbox{for }z\leq  -p_{1}(t),\\
&Q_z(\phi_{i,j;1}(z-p_{1}(t)),\phi_{i,j;2}(z+p_{1}(t)),\phi_{i,j;3}(z+p_{2}(t)))\geq \varepsilon_2 ,\ \ \mbox{for }  p_{1}(t)  \leq z\leq  -p_{2}(t),\\
&Q_w(\phi_{i,j;1}(z-p_{1}(t)),\phi_{i,j;2}(z+p_{1}(t)),\phi_{i,j;3}(z+p_{2}(t)))\geq \varepsilon_3 , \ \ \mbox{for }    z\geq  -p_{1}(t).\\
\end{align*}
\end{lemma}                               
From \eqref{2.9} and Lemma \eqref{2.10}, we see that
\begin{align*}
&\mathcal{F}(1,z,w)=f_{i,j}(Q(1,z,w))-Q_1'(1,z,w)f_{i,j}(1)-Q_2'(1,z,w)f_{i,j}(z)-Q_3'(1,z,w)f_{i,j}(w)=0,\\
&\mathcal{F}(y,0,w)=f_{i,j}(Q(y,0,w))-Q_1'(y,0,w)f_{i,j}(y)-Q_2'(y,0,w)f_{i,j}(0)-Q_3'(y,0,w)f_{i,j}(w)=0,\\
&\mathcal{F}(y,a,w)=f_{i,j}(Q(y,a,w))-Q_1'(y,a,w)f_{i,j}(y)-Q_2'(y,a,w)f_{i,j}(a)-Q_3'(y,a,w)f_{i,j}(w)=0,\\
&\mathcal{F}(y,z,1)=f_{i,j}(Q(y,z,1))-Q_1'(y,z,1)f_{i,j}(y)-Q_2'(y,z,1)f_{i,j}(z)-Q_3'(y,z,1)f_{i,j}(1)=0.\\
\end{align*}
Which implies that there is a smooth function $\mathcal{F}_1$ satisfying
\begin{align*}
\mathcal{F}(y,z,w)=(1-y)z(a-z)(1-w)\mathcal{F}_1(y,z,w).
\end{align*}
Since $Q(0,z,a)=z$ and $Q_2'(0,z,a)=1$, we have
\begin{align*}
\mathcal{F}(0,z,a)=f_{i,j}(Q(0,z,a))-Q_2'(0,z,a)f_{i,j}(z)=0,
\end{align*}
which implies $\mathcal{F}_1(0,z,a)=0$. Applying the mean value theorem to $\mathcal{F}_1$ yields
\begin{align*}
 \mathcal{F}_1(y,z,w)=&\int_{0}^1\mathcal{F}_{11}(\theta y,z,\theta w+(1-\theta)a)d\theta\cdot y \\+&\int_{0}^1\mathcal{F}_{13}(\theta y,z,\theta w+(1-\theta)a)d\theta\cdot (w-a)
\end{align*}
Thus we obtain
\begin{equation}
\left\{ {\begin{array}{ll}
\mathcal{F}(y,z,w)=(1-y)z[y\mathcal{F}_{11}(y,z,w)+(w-a)y\mathcal{F}_{12}(y,z,w)]\mbox \medskip \\
\mathcal{F}(y,z,w)=(1-w)(a-z)[y\mathcal{F}_{21}(y,z,w)+(w-a)y\mathcal{F}_{22}(y,z,w)]\mbox \medskip \\
\end{array}} \right.
\label{2.26}
\end{equation}
From the above discussion, we can easily obtain that there exists a positive constant $C$ such that
\begin{align*}
&|R_l(\phi_{i,j;1}, \phi_{i,j;2},\phi_{i,j;3})|, |\mathcal{F}_{m,n}(\phi_{i,j;1},\phi_{i,j;2},\phi_{i,j;3})|, |Q_{yy}(\phi_{i,j;1},\phi_{i,j;2},\phi_{i,j;3})|\\
&|Q_{zz}(\phi_{i,j;1},\phi_{i,j;2},\phi_{i,j;3})|, |Q_{ww}(\phi_{i,j;1},\phi_{i,j;2},\phi_{i,j;3})|, |Q_{yz}(\phi_{i,j;1},\phi_{i,j;2},\phi_{i,j;3})|\\
&|Q_{zw}(\phi_{i,j;1},\phi_{i,j;2},\phi_{i,j;3})|, |Q_{yw}(\phi_{i,j;1},\phi_{i,j;2},\phi_{i,j;3})|\leq C\\
\end{align*}
for $z\in R$,  $p_{1}(t)$,  $p_{2}(t)<-\delta$, $l=1,\cdot\cdot\cdot,16$, $m ,n=1,2.$\\
Set $\bar{c}:=(c_{1}+c_{2})/2$, $s_{1}:=(c_{2}-c_1)/2 $ and $s_{2}:=c_3-\bar{c}$. Define the functions $\overline{U}_{i,j}(t)$ and $\underline{U}_{i,j}(t)$ by
\begin{align*}
\overline{U}_{i,j}(t)=Q\big(&\phi_{i,j;1}(i cos\theta +j sin\theta+\bar{c}t-p_{1}(t)),\phi_{i,j;2}(i cos\theta
+jsin\theta+\bar{c}t+p_{1}(t)),\\
&\phi_{i,j;3}(i cos\theta +j sin\theta+\bar{c}t+p_{2}(t))\big),\\
\underline{U}_{i,j}(t)=Q\big(&\phi_{i,j;1}(i cos\theta +j sin\theta+\bar{c}t-r_{1}(t)),\phi_{i,j;2}(i cos\theta
+jsin\theta+\bar{c}t+r_{1}(t)),\\
&\phi_{i,j;3}(i cos\theta +j sin\theta+\bar{c}t+r_{2}(t))\big)            
\end{align*}
constitute a pair of super- and subsolution of \eqref{eq1.1} for $t\leq0$. Moreover, there exists a positive constant $C$ such that
\begin{align}
\overline{U}_{i,j}(t) >\underline{U}_{i,j}(t)\ \ \mbox{and}\ \
\sup_{i,j\in\mathbb{Z}}(\overline{U}_{i,j}(t)-\underline{U}_{i,j}(t))
\leq Ce^{ \kappa s_{1}t},~~\mbox{for }t\leq0,\label{2.27}
\end{align}

\paragraph{Proof.} We only prove that $\overline{U}_{i,j}(t)$ is a super-solution, since the other case can be discussed similarly. For convenience, we denote
\begin{align*}
&z=i cos\theta +j sin\theta+\bar{c}t\\
&z_1=i cos\theta +j sin\theta+\bar{c}t-p_1(t)=z-p_1(t)\\
&z_2=i cos\theta +j sin\theta+\bar{c}t+p_1(t)=z+p_1(t)\\
&z_3=i cos\theta +j sin\theta+\bar{c}t+p_2(t)=z+p_2(t)
\end{align*}
\begin{align*}
F(v)(i,j,t):=v_{i,j}'(t)-\sum_{k_1}\sum_{k_2}J(k_1,k_2)v_{i-k_1,j-k_2}(t)-v_{i,j}(t)-f_{i,j}(v_{i,j}(t)).          
\end{align*}
Direct computations show that
\begin{align*}
&F(\overline{U})(i,j,t)\nonumber\\
=& (\bar{c}-p_1')Q_1' \phi_{i,j;1}' +(\bar{c}+p_2')Q_2' \phi_{i,j;2}' +(\bar{c}+p_2')Q_3' \phi_{i,j;3}'-
\sum_{k_1}\sum_{k_2}J(k_1,k_2)\overline{U}_{i-k_1,j-k_2}(t)\\
 &+ \overline{U}_{i,j}(t)- f_{i,j}(\overline{U}_{i,j}(t))\nonumber\\
=& (\bar{c}-p_1'-c_1)Q_1' \phi_{i,j;1}' +(\bar{c}+p_2'-c_2)Q_2' \phi_{i,j;2}' +(\bar{c}+p_2'-c_3)Q_3'
\phi_{i,j;3}'\\
&-\sum_{k_1}\sum_{k_2}J(k_1,k_2)\overline{U}_{i-k_1,j-k_2}(t)+ \overline{U}_{i,j}(t)- f_{i,j}(\overline{U}_{i,j}(t))+c_1Q_1' \phi_{i,j;1}'+ c_2Q_2' \phi_{i,j;2}'\\
&+ c_3Q_3' \phi_{i,j;3}'\nonumber\\
=& (s_1-p_1')Q_1' \phi_{i,j,1}' +(-s_1+p_1')Q_2' \phi_{i,j,2}' +(-s_2+p_2')Q_3'\phi_{i,j,3}'\\
&-\sum_{k_1}\sum_{k_2}J(k_1,k_2)\overline{U}_{i-k_1,j-k_2}(t)+ \overline{U}_{i,j}(t)- f_{i,j}(\overline{U}_{i,j}(t))+c_1Q_1' \phi_{i,j;1}'+ c_2Q_2' \phi_{i,j;2}'\\
&+ c_3Q_3' \phi_{i,j;3}'\nonumber\\
=& -Le^{\kappa p_{1}}Q_1'\phi_{i,j;1}' +Le^{\kappa p_{1}}Q_2'\phi_{i,j;2}'+Le^{\kappa p_{2}}Q_3'
\phi_{i,j;3}'-\sum_{k_1}\sum_{k_2}J(k_1,k_2)\overline{U}_{i-k_1,j-k_2}(t)\\
 &+ \overline{U}_{i,j}(t)- f_{i,j}(\overline{U}_{i,j}(t))+c_1Q_1' \phi_{i,j;1}'+ c_2Q_2' \phi_{i,j;2}'+ c_3Q_3' \phi_{i,j;3}'\nonumber\\
=& -Le^{\kappa p_{1}}Q_1'\phi_{i,j;1}' +Le^{\kappa p_{1}}Q_2'\phi_{i,j;2}'+Le^{\kappa p_{2}}Q_3'
\phi_{i,j;3}'-\sum_{k_1}\sum_{k_2}J(k_1,k_2)\overline{U}_{i-k_1,j-k_2}(t)\\
& + \overline{U}_{i,j}(t)-f_{i,j}(\overline{U}_{i,j}(t))+Q_1'[\sum_{k_1}\sum_{k_2}J(k_1,k_2)\phi_{i-k_1,j-k_2;1}(t)-\phi_{i,j;1}(t)+f_{i,j}(\phi_{i,j;1}(t))] \nonumber\\
 &+Q_2'[\sum_{k_1}\sum_{k_2}J(k_1,k_2)\phi_{i-k_1,j-k_2;2}(t)-\phi_{i,j;2}(t)+f_{i,j}(\phi_{i,j;2}(t))]\nonumber \\
 &+Q_3'[\sum_{k_1}\sum_{k_2}J(k_1,k_3)\phi_{i-k_1,j-k_2;3}(t)-\phi_{i,j;3}(t)+f_{i,j}(\phi_{i,j;3}(t))]\nonumber\\
=& -Le^{\kappa p_{1}}Q_1'\phi_{i,j;1}' +Le^{\kappa p_{1}}Q_2'\phi_{i,j;2}'+Le^{\kappa p_{2}}Q_3' \phi_{i,j;3}'-
\mathcal{F}(\phi_{i,j;1},\phi_{i,j;2},\phi_{i,j;3})\\
&-\mathcal{H}(\phi_{i,j;1},\phi_{i,j;2},\phi_{i,j;3})\nonumber\\
\label{F}
\end{align*}
where $Q:=Q(\phi_{i,j;1}(z_1), \phi_{i,j;2}(z_2),\phi_{i,j;3}(z_3))$, $\phi_{i,j;1}:=\phi_{i,j;1}(z_1)$,
$\phi_{i,j;2}:=\phi_{i,j;2}(z_2)$, $\phi_{i,j;3}:=\phi_{i,j;3}(z_3)$, $\phi_{i-k_1,j-k_2;1}=\phi_{i-k_1,j-k_2;1}(z_1-k_1cos\theta-k_2sin\theta)$,
$\phi_{i-k_1,j-k_2;2}=\phi_{i-k_1,j-k_2;2}(z_2-k_1cos\theta-k_2sin\theta)$, $\phi_{i-k_1,j-k_2;3}=\phi_{i-k_1,j-k_2;3}(z_3-k_1cos\theta-k_2sin\theta)$ and
\begin{align*}
\mathcal{F}(\phi_{i,j;1},\phi_{i,j;2},\phi_{i,j;3}):=&f_{i,j}(\overline{U}_{i,j}(t))-Q_1'f_{i,j}(\phi_{i,j;1})-Q_2'f_{i,j}(\phi_{i,j;2})-Q_3'f_{i,j}(\phi_{i,j;3})\\
\mathcal{H}(\phi_{i,j;1},\phi_{i,j;2},\phi_{i,j;3}):=&\sum_{k_1}\sum_{k_2}J(k_1,k_2)\overline{U}_{i-k_1,j-k_2}(t)-\overline{U}_{i,j}(t)\\
-&Q_1'[\sum_{k_1}\sum_{k_2}J(k_1,k_2)\phi_{i-k_1,j-k_2;1}(t)-\phi_{i,j;1}(t)] \nonumber\\
-&Q_2'[\sum_{k_1}\sum_{k_2}J(k_1,k_2)\phi_{i-k_1,j-k_2;2}(t)-\phi_{i,j;2}(t)]\nonumber \\
-&Q_3'[\sum_{k_1}\sum_{k_2}J(k_1,k_3)\phi_{i-k_1,j-k_2;3}(t)-\phi_{i,j;3}(t)]\nonumber\\
\end{align*}
Denote \begin{align*}
&\mathcal{A}(\phi_{i,j;1}(z-p_1(t)),\phi_{i,j;2}(z+p_1(t)),\phi_{i,j;3}(z+p_2(t)))\\
&:=-Q_y(\phi_{i,j;1},\phi_{i,j;2},\phi_{i,j;3})\phi_{i,j;1}'+Q_z(\phi_{i,j;1},
\phi_{i,j;2},\phi_{i,j;3})\phi_{i,j;2}'+Q_w(\phi_{i,j;1}, \phi_{i,j;2},\phi_{i,j;3})\phi_{i,j;3}'\\
\end{align*}
\begin{lemma}\label{lem2.4}
There exists a large $\sigma>0$. The following statements hold that
\begin{equation}
\mathcal{A}(\phi_{i,j;1}, \phi_{i,j;2},\phi_{i,j;3})>0,\;\text{ for}\; z \in \mathbb{R},\;
p_1,p_2\le-\sigma.
\label{eq2.29}
\end{equation}
Moreover, if $p_1 ,p_2 \leq -\sigma $, then the following assertions hold:
\begin{align}
&\mathcal{A}(\phi_{i,j;1}, \phi_{i,j;2},\phi_{i,j;3})\geq\frac{1}{2}Q_y|\phi_{i,j;1}'(z-p_{1}(t))| \; for\;z\leq p_{1}(t),
\label{yin1.1}\\
&\mathcal{A}(\phi_{i,j;1}, \phi_{i,j;2},\phi_{i,j;3})\geq\frac{1}{2}[Q_y|\phi_{i,j;1}'(z-p_{1}(t))|+Q_z|\phi_{i,j;2}'(z+p_{1}(t))|]\nonumber\\
&\; for\;  p_{1}(t)  \leq z\leq -p_{1}(t),\label{yin1.2}\\
&\mathcal{A}(\phi_{i,j;1}, \phi_{i,j;2},\phi_{i,j;3})\geq\frac{1}{2}[Q_z|\phi_{i,j;2}'(z+p_{1}(t))|+Q_w|\phi_{i,j;3}'(z+p_{1}(t))|]\nonumber\\
&\; for \;  -p_{1}(t)  \leq z\leq -p_{2}(t),\label{yin1.3}\\
&\mathcal{A}(\phi_{i,j;1}, \phi_{i,j;2},\phi_{i,j;3})\geq \frac{1}{2}Q_w|\phi_{i,j;3}'(z+p_{1}(t))| \; for \; z\geq -p_{2}(t).
\label{yin1.4}
\end{align}
\end{lemma}
It follows from that
\begin{align*}
F(\overline{U})(i,j,t)\geq\mathcal{A}(\phi_{i,j;1}, \phi_{i,j;2},\phi_{i,j;3})[Le^{\kappa p_{1}}-\frac{
\mathcal{F}(\phi_{i,j;1},\phi_{i,j;2},\phi_{i,j;3})+\mathcal{H}(\phi_{i,j;1},\phi_{i,j;2},\phi_{i,j;3})}{\mathcal{A}(\phi_{i,j;1}, \phi_{i,j;2},\phi_{i,j;3})}].
\end{align*}
By direct computations, we obtain that
\begin{align*}
 &\sum_{k_1}\sum_{k_2}J(k_1,k_2)\overline{U}_{i-k_1,j-k_2}(t)-\overline{U}_{i,j}(t)\nonumber\\
 = &\sum_{k_1}\sum_{k_2}J(k_1,k_2)[Q(\phi_{i-k_1,j-k_2;1}(t),\phi_{i-k_1,j-k_2;2}(t),\phi_{i-k_1,j-k_2;3}(t))\\
 &-Q(\phi_{i,j;1}(t),\phi_{i-k_1,j-k_2;2}(t),\phi_{i-k_1,j-k_2;3}(t))
 +Q(\phi_{i,j;1}(t),\phi_{i-k_1,j-k_2;2}(t),\phi_{i-k_1,j-k_2;3}(t))\\
 &-Q(\phi_{i,j;1}(t),\phi_{i,j;2}(t),\phi_{i-k_1,j-k_2;3}(t))+Q(\phi_{i,j;1}(t),\phi_{i,j;2}(t),\phi_{i-k_1,j-k_2;3}(t)\\
 &-Q(\phi_{i,j;1}(t),\phi_{i,j;2}(t),\phi_{i,j;3}(t)]\nonumber\\
=&\sum_{k_1}\sum_{k_2}J(k_1,k_2)
[Q_1'(\eta^{(1)}_{i-k_1,j-k_2},\phi_{i-k_1,j-k_2;2}(t),\phi_{i-k_1,j-k_2;3}(t))[\phi_{i-k_1,j-k_2;1}(t)-\phi_{i,j;1}(t)]\\
&+[Q_2'(\phi_{i,j;1}(t),\eta^{(2)}_{i-k_1,j-k_2},\phi_{i-k_1,j-k_2;3}(t))[\phi_{i-k_1,j-k_2;2}(t)-\phi_{i,j;2}(t)]\\
&+[Q_3'(\phi_{i,j;1}(t),\phi_{i,j;2}(t),\eta^{(3)}_{i-k_1,j-k_2})[\phi_{i-k_1,j-k_2;3}(t)-\phi_{i,j;3}(t)]\\
 =&\sum_{k_1}\sum_{k_2}J(k_1,k_2)
[Q_{11}''(\eta^{(4)}_{i-k_1,j-k_2},\phi_{i-k_1,j-k_2;2}(t),\phi_{i-k_1,j-k_2;3}(t))(\eta^{(1)}_{i-k_1,j-k_2}-\phi_{i,j;1}(t)) \\
&+Q_{12}''(\phi_{i,j;1}(t),\eta^{(5)}_{i-k_1,j-k_2},\phi_{i-k_1,j-k_2;3}(t))(\phi_{i-k_1,j-k_2;2}(t)-\phi_{i,j;2}(t)) \\
&+Q_{13}''(\phi_{i,j;1}(t),\phi_{i,j;2}(t),\eta^{(5)}_{i-k_1,j-k_2})(\phi_{i-k_1,j-k_2;3}(t)-\phi_{i,j;3}(t))]
[\phi_{i-k_1,j-k_2;1}(t)-\phi_{i,j;1}(t)]\\
&+\sum_{k_1}\sum_{k_2}J(k_1,k_2)[Q_{22}''(\phi_{i,j;1}(t),\eta^{(7)}_{i-k_1,j-k_2},\phi_{i-k_1,j-k_2;3}(t))(\eta^{(2)}_{i-k_1,j-k_2}-\phi_{i,j;2}(t))
\\
&+Q_{23}''(\phi_{i,j;1}(t),\phi_{i,j;2}(t),\eta^{(8)}_{i-k_1,j-k_2})(\phi_{i-k_1,j-k_2;3}(t)-\phi_{i,j;3}(t))]
[\phi_{i-k_1,j-k_2;2}(t)-\phi_{i,j;2}(t)]\\
&+\sum_{k_1}\sum_{k_2}J(k_1,k_2)[Q_{33}''(\phi_{i,j;1}(t),\phi_{i,j;2}(t),\eta^{(9)}_{i-k_1,j-k_2})(\eta^{(3)}_{i-k_1,j-k_2}-\phi_{i,j;3}(t)]
[\phi_{i-k_1,j-k_2;3}(t)\\
&-\phi_{i,j;3}(t)]\\
\end{align*}

\begin{align*}
=&\sum_{k_1}\sum_{k_2}J(k_1,k_2)
[Q_{11}''(\eta^{(4)}_{i-k_1,j-k_2},\phi_{i-k_1,j-k_2;2}(t),\phi_{i-k_1,j-k_2;3}(t))\theta^{(1)}_{i-k_1,j-k_2}(\phi_{i-k_1,j-k_2;1}(t)\\
&-\phi_{i,j;1}(t))^2+\sum_{k_1}\sum_{k_2}J(k_1,k_2)Q_{12}''(\phi_{i,j;1}(t),\eta^{(5)}_{i-k_1,j-k_2},\phi_{i-k_1,j-k_2;3}(t))(\phi_{i-k_1,j-k_2;2}(t)\\
&-\phi_{i,j;2}(t))(\phi_{i-k_1,j-k_2;1}(t)-\phi_{i,j;1}(t))+\sum_{k_1}\sum_{k_2}J(k_1,k_2)Q_{13}''(\phi_{i,j;1}(t),\phi_{i,j;2}(t),\eta^{(5)}_{i-k_1,j-k_2})\\
&(\phi_{i-k_1,j-k_2;3}(t)-\phi_{i,j;3}(t))(\phi_{i-k_1,j-k_2;1}(t)-\phi_{i,j;1}(t))+\sum_{k_1}\sum_{k_2}J(k_1,k_2)[Q_{22}''(\phi_{i,j;1}(t),\\
&\eta^{(7)}_{i-k_1,j-k_2},\phi_{i-k_1,j-k_2;3}(t))\theta^{(2)}_{i-k_1,j-k_2}(\phi_{i-k_1,j-k_2;2}(t)-\phi_{i,j;2}(t))^2+\sum_{k_1}\sum_{k_2}J(k_1,k_2)\\
&Q_{23}''(\phi_{i,j;1}(t),\phi_{i,j;2}(t),\eta^{(8)}_{i-k_1,j-k_2})(\phi_{i-k_1,j-k_2;3}(t)-\phi_{i,j;3}(t))
(\phi_{i-k_1,j-k_2;2}(t)-\phi_{i,j,2}(t))\\
&+\sum_{k_1}\sum_{k_2}J(k_1,k_2)[Q_{33}''(\phi_{i,j;1}(t),\phi_{i,j;2}(t),\eta^{(9)}_{i-k_1,j-k_2})\theta^{(3)}_{i-k_1,j-k_2}(\phi_{i-k_1,j-k_2;3}(t)-\phi_{i,j;3}(t))^2\\
\end{align*}

where
\begin{align*}
 &\eta^{(1)}_{i-k_1,j-k_2}:=\theta^{(1)}_{i-k_1,j-k_2}\phi_{i-k_1,j-k_2;1}(t)+(1-\theta^{(1)}_{i-k_1,j-k_2})\phi_{i,j;1}(t),\\
 &\eta^{(2)}_{i-k_1,j-k_2}:=\theta^{(2)}_{i-k_1,j-k_2}\phi_{i-k_1,j-k_2;2}(t)+(1-\theta^{(2)}_{i-k_1,j-k_2})\phi_{i,j;2}(t),\\
& \eta^{(3)}_{i-k_1,j-k_2}:=\theta^{(3)}_{i-k_1,j-k_2}\phi_{i-k_1,j-k_2;3}(t)+(1-\theta^{(3)}_{i-k_1,j-k_2})\phi_{i,j;3}(t),\\
&\eta^{(4)}_{i-k_1,j-k_2}:=\theta^{(4)}_{i-k_1,j-k_2}\eta^{(1)}_{i-k_1,j-k_2}+(1-\theta^{(4)}_{i-k_1,j-k_2})\phi_{i,j;1}(t),\\
&\eta^{(5)}_{i-k_1,j-k_2}:=\theta^{(5)}_{i-k_1,j-k_2}\eta^{(2)}_{i-k_1,j-k_2}+(1-\theta^{(5)}_{i-k_1,j-k_2})\phi_{i,j;2}(t),\\
&\eta^{(6)}_{i-k_1,j-k_2}:=\theta^{(6)}_{i-k_1,j-k_2}\eta^{(3)}_{i-k_1,j-k_2}+(1-\theta^{(6)}_{i-k_1,j-k_2})\phi_{i,j;3}(t),\\
&\eta^{(1)}_{i-k_1,j-k_2}-\phi_{i,j;1}(t):=\theta^{(1)}_{i-k_1,j-k_2}(\phi_{i-k_1,j-k_2;1}(t)-\phi_{i,j;1}(t))\\
&\eta^{(2)}_{i-k_1,j-k_2}-\phi_{i,j;2}(t):=\theta^{(2)}_{i-k_1,j-k_2}(\phi_{i-k_1,j-k_2;2}(t)-\phi_{i,j;2}(t))\\
&\eta^{(3)}_{i-k_1,j-k_2}-\phi_{i,j;3}(t):=\theta^{(3)}_{i-k_1,j-k_2}(\phi_{i-k_1,j-k_2;3}(t)-\phi_{i,j;3}(t))\\
\end{align*}
There is a positive constant $M$ such that
\begin{align*}
&|\frac{\mathcal{H}(\phi_{i,j;1},\phi_{i,j;2},\phi_{i,j;3})}{ \mathcal{A}(\phi_{i,j;1},\phi_{i,j;2},\phi_{i,j;3})}|\\
&\leq\left\{ {\begin{array}{ll}
L_1(e^{\eta_{1}p_{1}(t)} + e^{\eta_{1}p_{2}(t)}), \; for\; z\leq0 \mbox \medskip \\
L_2(e^{\eta_{2}p_{1}(t)} + e^{\eta_{1}(p_{2}(t)-p_{1}(t))/2}), \; for\; 0\leq z\leq(-p_{2}-p_{1}(t))/2 \mbox \medskip \\
L_3(e^{\eta_{2}p_{1}(t)} + e^{\eta_{2}(p_{2}(t)-p_{1}(t))/2}), \; for\; (-p_{2}-p_{1}(t))/2 \leq z\leq-p_{2}\mbox \medskip \\
 \end{array}} \right.
\end{align*}
When $z\leq p_1(t)$, based on Lemma \ref{lem2.1}, Lemma \ref{lem2.2}, Lemma \ref{lem2.3} and \eqref{yin1.1}, we have $z_1=z-p_{1}(t)\leq0$,  $z_2=z+p_{1}(t)\leq0$
and $ z_3=z+p_{2}(t)\leq0$. Then it follows that
\begin{align*}
&\frac{   \mathcal{H}(\phi_{i,j;1},\phi_{i,j;2},\phi_{i,j;3})}{ \mathcal{A}(\phi_{i,j;1},\phi_{i,j;2},\phi_{i,j;3})}\\
\leq&\sum_{k_1}\sum_{k_2}J(k_1,k_2)  \frac{  \tilde{C}_1
|\phi_{i-k_1,j-k_2;2}|(\phi_{i-k_1,j-k_2;1}-\phi_{i,j;1})^2}{\frac{1}{2}Q_1'|\phi_{i,j;1}'|}\\
&+\sum_{k_1}\sum_{k_2}J(k_1,k_2) \frac{ \tilde{C}_2
(\phi_{i-k_1,j-k_2;2}-\phi_{i,j;2})(\phi_{i-k_1,j-k_2;1}-\phi_{i,j;1})}{\frac{1}{2}Q_1'|\phi_{i,j;1}'|}\\
&+\sum_{k_1}\sum_{k_2}J(k_1,k_2) \frac{\tilde{C}_3
(\phi_{i-k_1,j-k_2;3}-\phi_{i,j,3})(\phi_{i-k_1,j-k_2;1}-\phi_{i,j;1})}{\frac{1}{2}Q_1'|\phi_{i,j;1}'|}\\
&+\sum_{k_1}\sum_{k_2}J(k_1,k_2)\frac{ \tilde{C}_4
(1-\phi_{i,j;1})(\phi_{i-k_1,j-k_2;2}-\phi_{i,j;2})^2}{\frac{1}{2}Q_1'|\phi_{i,j;1}'|}\\
&+\sum_{k_1}\sum_{k_2}J(k_1,k_2)\frac{ \tilde{C}_5
(1-\phi_{i,j;1})(\phi_{i-k_1,j-k_2;3}-\phi_{i,j;3})(\phi_{i-k_1,j-k_2;2}-\phi_{i,j;2})}{\frac{1}{2}Q_1'|\phi_{i,j;1}'|}\\
&+\sum_{k_1}\sum_{k_2}J(k_1,k_2)\frac{ \tilde{C}_6
(1-\phi_{i,j;1})(\phi_{i-k_1,j-k_2;3}-\phi_{i,j;3})^2}{\frac{1}{2}Q_1'|\phi_{i,j;1}'|}\\
\leq&\frac{\tilde{C}_1C_2e^{z_2-k_1 cos\theta -k_2 sin\theta}[ C_2e^{\eta_1(z_1-k_1 cos\theta -k_2 sin\theta)}-C_1e^{\eta_1z_1}
]^2}{\epsilon_1\rho C_1 e^{\eta_1z_1}}\\
&+\frac{\tilde{C}_2[C_2e^{\eta_1(z_2-k_1cos\theta -k_2 sin\theta)}-C_1 e^{\eta_1z_2}][C_2e^{\eta_1(z_1-k_1cos\theta-k_2
sin\theta)}-C_1e^{\eta_1z_1}]}{\epsilon_1\rho C_1 e^{\eta_1z_1}  }\\
&+\frac{\tilde{C}_3[ C_2e^{\eta_1(z_3-k_1 cos\theta -k_2 sin\theta)}-C_1e^{\eta_1z_3}][C_2e^{\eta_1(z_1-k_1cos\theta-k_2
sin\theta)}-C_1e^{\eta_1z_1}]}{\epsilon_1\rho C_1 e^{\eta_1z_1}  }\\
&+\frac{\tilde{C}_4C_2e^{\eta_1 z_1}[ C_2e^{\eta_1(z_2-k_1 cos\theta -k_2 sin\theta)}-C_1 e^{\eta_1z_2}]^2}{\epsilon_1\rho C_1e^{\eta_1z_1}}\\
&+\frac{\tilde{C}_5C_1e^{\eta_1z_1}[C_2e^{\eta_1(z_3-k_1cos\theta -k_2 sin\theta)}-C_1e^{\eta_1z_3}][C_2e^{\eta_1(z_2-k_1
cos\theta -k_2 sin\theta)}-C_1e^{\eta_1z_2}]}{\epsilon_1\rho C_1 e^{\eta_1z_1}  }\\
&+\frac{\tilde{C}_6C_1e^{\eta_1z_1}[C_2e^{\eta_1(z_3-k_1 cos\theta -k_2 sin\theta)}-C_1 e^{\eta_1z_3} ]^2}{\epsilon_1\rho C_1
e^{\eta_1z_1}}\\
\leq&\tilde{L}_1(e^{\eta_{1}(z+p_{1}(t))}+e^{\eta_{1}(z+p_{2}(t))})\leq\tilde{L}_1(e^{\eta_{1}p_{1}(t)}+e^{\eta_{1}p_{2}(t)})\\
\end{align*}

For $p_1(t)\leq z\leq0$, based on Lemma \ref{lem2.1}, Lemma \ref{lem2.2}, Lemma \ref{lem2.3} and \eqref{yin1.2},
we have $z_1=z-p_{1}(t)\geq0$, $z_2=z+p_{1}(t)\leq0$ and $ z_3=z+p_{2}(t)\leq0$. Then it follows that
\begin{align*}
&|\frac{\mathcal{H}(\phi_{i,j;1},\phi_{i,j;2},\phi_{i,j;3})}{ \mathcal{A}(\phi_{i,j;1},\phi_{i,j;2},\phi_{i,j;3})}|\\
\leq&\frac{\hat {C}_1|\phi_{i-k_1,j-k_2;2}|(\phi_{i-k_1,j-k_2;1}-\phi_{i,j;1})^2}{\frac{1}{2}Q_1'|\phi_{i,j;1}'|}\\
&+\frac{\hat{C}_2(\phi_{i-k_1,j-k_2;2}-\phi_{i,j;2})(\phi_{i-k_1,j-k_2;1}-\phi_{i,j;1})}{\frac{1}{2}Q_1'|\phi_{i,j;1}'|}\\
&+\frac{\hat{C}_3(\phi_{i-k_1,j-k_2;3}-\phi_{i,j;3})(\phi_{i-k_1,j-k_2;1}-\phi_{i,j;1})}{\frac{1}{2}Q_1'|\phi_{i,j;1}'|}\\
&+\frac{\hat{C}_4(1-\phi_{i-k_1,j-k_2;1})(\phi_{i-k_1,j-k_2;2}-\phi_{i,j;2})^2}{\frac{1}{2}Q_2'|\phi_{i,j;2}'|}\\
&+\frac{\hat{C}_5(1-\phi_{i-k_1,j-k_2;1})(\phi_{i-k_1,j-k_2;3}-\phi_{i,j;3})(\phi_{j+1,2},-\phi_{j,2})}{\frac{1}{2}Q_2'|\phi_{i,j;2}'|}\\
&+\frac{\hat {C}_6 (1-\phi_{i-k_1,j-k_2;1})(\phi_{i-k_1,j-k_2;3}-\phi_{i,j;3})^2}{\frac{1}{2}Q_2'|\phi_{i,j;2}'|}\\
\leq&\frac{\hat {C}_1C_2e^{\eta_1 (z_2-k_1 cos\theta -k_2 sin\theta)}[ C_2e^{\eta_1(z_1-k_1 cos\theta -k_2 sin\theta)}-C_1
e^{\eta_1z_1}]^2}{\epsilon_1\rho C_1e^{-\eta_2z_1}}\\
&+\frac{\hat{C}_2[C_2e^{\eta_1(z_2-k_1 cos\theta -k_2 sin\theta)}-C_1e^{\eta_1z_2}][C_2e^{\eta_1(z_1-k_1 cos\theta-k_2
sin\theta)}-C_1e^{\eta_1z_1}]}{\epsilon_1\rho C_1 e^{-\eta_2z_1}}\\
&+\frac{\hat{C}_3[C_2e^{\eta_1(z_3-k_1 cos\theta -k_2 sin\theta)}-C_1e^{\eta_1z_3}][C_2e^{\eta_1(z_1-k_1 cos\theta-k_2
sin\theta)}-C_1e^{\eta_1z_1}]}{\epsilon_1\rho C_1 e^{-\eta_2z_1}}\\
&+\frac{\hat {C}_4[C_2e^{\eta_1(z_2-k_1 cos\theta -k_2 sin\theta)}-C_1e^{\eta_1z_2} ]^2}{\epsilon_1\rho C_1 e^{\eta_1z_2}}\\
&+\frac{\hat {C}_5[\hat C_2e^{\eta_1(z_3-k_1 cos\theta -k_2 sin\theta)}-C_1 e^{\eta_1z_3}][C_2e^{\eta_1(z_2-k_1 cos\theta
-k_2 sin\theta)}-C_1e^{\eta_1z_2}]}{\epsilon_1\rho C_1 e^{\eta_1z_2}}\\
&+\frac{\hat {C}_6 C_1e^{\eta_1z_1}[ C_2e^{\eta_1(z_3-k_1 cos\theta -k_2 sin\theta)}-C_1 e^{\eta_1z_3} ]^2}{\epsilon_1\rho C_1
e^{\eta_1z_2}}\\
\leq&\tilde{L}_2(e^{\eta_{1}(z+p_{1}(t))}+e^{\eta_{1}(z+p_{2}(t))})\leq\tilde{L}_2(e^{\eta_{1}p_{1}(t)} +e^{\eta_{1}p_{2}(t)})\\
\end{align*}
Then, we can obtain that
 \begin{align*}
|\frac{\mathcal{H}(\phi_{i,j;1}, \phi_{i,j;2},\phi_{i,j;3})}{ \mathcal{A}(\phi_{i,j;1}, \phi_{i,j;2},\phi_{i,j;3})}|\leq
L_1(e^{\eta_{1}p_{1}(t)} + e^{\eta_{1}p_{2}(t)}), \; for\; z\leq0
\end{align*}
where $L_1=max\{\tilde{L}_1,\tilde{L}_2\}$.

For $0\leq z\leq-p_{1}(t)$, based on Lemma \ref{lem2.1}, Lemma \ref{lem2.2}, Lemma \ref{lem2.3} and \eqref{yin1.2},
we have $z_1=z-p_{1}(t)\geq0$, $z_2=z+p_{1}(t)\leq0$ and $ z_3=z+p_{2}(t)\leq0$. Then it follows that
\begin{align*}
&|\frac{\mathcal{H}(\phi_{i,j;1},\phi_{i,j;2},\phi_{i,j;3})}{\mathcal{A}(\phi_{i,j;1},\phi_{i,j;2},\phi_{i,j;3})}|\\
\leq&\frac{\check{C}_1(\phi_{i-k_1,j-k_2;1}-\phi_{i,j;1}{z_1})^2}{\frac{1}{2}Q_1'|\phi_{i,j;1}'|}\\
&+\frac{\check{C}_2(\phi_{i-k_1,j-k_2;2}-\phi_{i,j;2}{z_2})(\phi_{i-k_1,j-k_2;1}-\phi_{i,j;1}{z_1})}{\frac{1}{2}Q_2'|\phi_{i,j;2}'|}\\
&+\frac{\check{C}_3(\phi_{i-k_1,j-k_2;3}-\phi_{i,j;3}{z_3})(\phi_{i-k_1,j-k_2;1}-\phi_{i,j;1}{z_1})}{\frac{1}{2}Q_1'|\phi_{i,j;1}'|}\\
&+\frac{\check{C}_4(\phi_{i,j;1}R_6+(\phi_{i-k_1,j-k_2;3}-a)R_7)(\phi_{j+1;2},-\phi_{i-k_1,j-k_2;2})^2}{\frac{1}{2}Q_2'|\phi_{i,j;2}'|}\\
&+\frac{\check{C}_5(\phi_{i-k_1,j-k_2;3}-\phi_{i,j;3})(\phi_{i-k_1,j-k_2;2}-\phi_{i,j;2})}{\frac{1}{2}Q_2'|\phi_{i,j;2}'|}\\
&+\frac{\check{C}_6(\phi_{i,j;2}R_9)(\phi_{i-k_1,j-k_2;3}-\phi_{i,j;3})^2}{\frac{1}{2}Q_2'|\phi_{i,j;2}'|}\\
\leq&\frac{\check{C}_1[ C_2e^{-\eta_2(z_1-k_1 cos\theta -k_2 sin\theta)}-C_1 e^{-\eta_2z_1} ]^2}{\epsilon_1\rho C_1e^{-\eta_2z_1}}\\
&+\frac{\check{C}_2[C_2e^{\eta_1(z_2-k_1cos\theta-k_2sin\theta)}-C_1e^{\eta_1z_2}][C_2e^{-\eta_2(z_1-k_1cos\theta-k_2sin\theta)}-C_1e^{-\eta_2z_1}]}{\epsilon_1\rho C_1 e^{\eta_1z_2}}\\
&+\frac{\check{C}_3[C_2e^{\eta_1(z_3-k_1 cos\theta -k_2 sin\theta)}-C_1e^{\eta_1z_3}][C_2e^{-\eta_2(z_1-k_1 cos\theta -k_2
sin\theta)}-C_1e^{-\eta_2z_1}]}{\epsilon_1\rho C_1 e^{-\eta_2z_1}}\\
&+\frac{\check{C}_4[C_2e^{-\eta_2z_1}R_6+C_1 e^{\eta_1(z_3-k_1 cos\theta -k_2 sin\theta)}R_7][C_2e^{\eta_1(z_2-k_1 cos\theta
-k_2 sin\theta)}-C_1e^{\eta_1z_2}]^2}{\epsilon_1\rho C_1e^{\eta_1z_2}}\\
&+\frac{\check{C}_5[ C_2e^{\eta_1(z_3-k_1 cos\theta -k_2 sin\theta)}-C_1 e^{\eta_1z_3}][C_2e^{\eta_1(z_2-k_1 cos\theta -k_2
sin\theta)}-C_1e^{\eta_1z_2}]}{\epsilon_1\rho C_1 e^{\eta_1z_2}}\\
&+\frac{\check{C}_6 R_9 C_2e^{\eta_1z_2}[ C_2e^{\eta_1(z_3-k_1 cos\theta -k_2 sin\theta)}-C_1 e^{\eta_1z_3}]}{\epsilon_1\rho
C_1e^{\eta_1z_2}}\\
\leq&\tilde{L}_3(e^{-\eta_{2}(z-p_{1}(t))}+e^{\eta_{1}(z+p_{2}(t))})\leq\tilde{L}_3(e^{\eta_{2}p_{1}(t)}+
e^{\eta_{1}(p_{2}(t)-p_{1}(t))/2})\\
\end{align*}

For $-p_{1}(t)\leq z\leq(-p_{2}-p_{1}(t))/2$, based on Lemma \ref{lem2.1}, Lemma \ref{lem2.2}, Lemma \ref{lem2.3} and \eqref{yin1.3},
we have $z_1=z-p_{1}(t)\geq0$, $z_2=z+p_{1}(t)\geq0$ and $z_3=z+p_{2}(t)\leq0$. Then it follows that
\begin{align*}
&|\frac{\mathcal{H}(\phi_{i,j;1},\phi_{i,j;2},\phi_{i,j;3})}{ \mathcal{A}(\phi_{i,j;1},\phi_{i,j;2},\phi_{i,j;3})}|\\
\leq&\frac{\bar{C}_1 (a-\phi_{i-k_1,j-k_2;2})(\phi_{i-k_1,j-k_2;1}-\phi_{i,j;1})^2}{\frac{1}{2}Q_2'|\phi_{i,j;2}'|}\\
&+\frac{\bar{C}_2(\phi_{i-k_1,j-k_2;2}-\phi_{i,j;2})(\phi_{i-k_1,j-k_2;1}-\phi_{i,j;1})}{\frac{1}{2}Q_2'|\phi_{i,j;2}'|}\\
&+\frac{\bar{C}_3(\phi_{i-k_1,j-k_2;3}-\phi_{i,j;3})(\phi_{i-k_1,j-k_2;1}-\phi_{i,j;1})}{\frac{1}{2}Q_3'|\phi_{i,j;3}'|}\\
&+\frac{\bar{C}_4(\phi_{i,j;1}R_6+(\phi_{i-k_1,j-k_2;3}-a)R_7)(\phi_{i-k_1,j-k_2;2}-\phi_{i,j;2})^2}{\frac{1}{2}Q_2'|\phi_{i,j;2}'|}\\
&+\frac{\bar{C}_5(\phi_{i-k_1,j-k_2;3}-\phi_{i,j;3})(\phi_{i-k_1,j-k_2;2}-\phi_{i,j;2})}{\frac{1}{2}Q_2'|\phi_{i,j;2}'|}\\
&+\frac{\bar{C}_6(\phi_{i-k_1,j-k_2;3}-\phi_{i,j;3})^2}{\frac{1}{2}Q_3'|\phi_{i,j;3}'|}\\
\leq&\frac{\bar{C}_1C_2e^{-\eta_2(z_2-k_1 cos\theta -k_2 sin\theta)}[C_2e^{-\eta_2(z_1-k_1cos\theta-k_2sin\theta)}- C_1
e^{-\eta_2z_1}]^2}{\epsilon_1\rho C_1 e^{-\eta_2z_2}}\\
&+\frac{\bar{C}_2[C_2e^{-\eta_2(z_2-k_1 cos\theta -k_2 sin\theta)}-C_1e^{-\eta_2z_2}][C_2e^{-\eta_2(z_1-k_1cos\theta-k_2
sin\theta)}-C_1e^{-\eta_2z_1}]}{\epsilon_1\rho C_1 e^{-\eta_2z_2}}\\
&+\frac{\bar{C}_3[C_2e^{\eta_1(z_3-k_1 cos\theta -k_2 sin\theta)}-C_1e^{\eta_1z_3}][C_2e^{-\eta_2(z_1-k_1cos\theta-k_2
sin\theta)}-C_1e^{-\eta_2z_1}]}{\epsilon_1\rho C_1 e^{-\eta_1z_3}}\\
&+\frac{\bar{C}_4[C_2e^{-\eta_2z_1}R_6+C_1 e^{\eta_1(z_3-k_1 cos\theta-k_2 sin\theta)}R_7 ][ C_2e^{-\eta_2(z_2-k_1 cos\theta
-k_2sin\theta)}-C_1 e^{-\eta_2z_2}]}{\epsilon_1\rho C_1 e^{-\eta_2z_2}}\\
&+\frac{\bar{C}_5[C_2e^{\eta_1(z_3-k_1 cos\theta -k_2 sin\theta)}-C_1e^{\eta_1z_3}][C_2e^{-\eta_2(z_2-k_1 cos\theta -k_2
sin\theta)}-C_1e^{-\eta_2z_2}]}{\epsilon_1\rho C_1 e^{-\eta_2z_2}}\\
&+\frac{\bar{C}_6[C_2e^{\eta_1(z_3-k_1 cos\theta -k_2 sin\theta)}-C_1 e^{\eta_1z_3}]}{\epsilon_1\rho C_1 e^{-\eta_1z_3}}\\
\leq&\tilde{L}_4 (e^{-\eta_{2}(z-p_{1}(t))}+e^{\eta_{1}(z+p_{2}(t))})\leq\tilde{L}_4(e^{\eta_{2}p_{1}(t)}+e^{\eta_{1}(p_{2}(t)-p_{1}(t))/2})\\
\end{align*}
Then, we can obtain that
 \begin{align*}
|\frac{\mathcal{H}(\phi_{i,j;1}, \phi_{i,j;2},\phi_{i,j;3})}{ \mathcal{A}(\phi_{i,j;1}, \phi_{i,j;2},\phi_{i,j;3})}|\leq
L_2(e^{-\eta_{2}(z-p_{1}(t))} + e^{\eta_{1}(z+p_{2}(t))}),, \; for\; 0\leq z\leq(-p_{2}-p_{1}(t))/2,
\end{align*}
where $L_2=max\{\tilde{L}_3,\tilde{L}_4\}$.

For $(-p_{2}-p_{1}(t))/2\leq z\leq-p_{2} $, based on Lemma\ref{lem2.1}, Lemma \ref{lem2.2}, Lemma \ref{lem2.3} and \eqref{yin1.3},
we have $z_1=z-p_{1}(t)\geq0$, $z_2=z+p_{1}(t)\geq0$ and $z_3=z+p_{2}(t)\leq0$. Then it follows that
\begin{align*}
&|\frac{\mathcal{H}(\phi_{i,j;1},\phi_{i,j;2},\phi_{i,j;3})}{\mathcal{A}(\phi_{i,j;1},\phi_{i,j;2},\phi_{i,j;3})}|\\
\leq&\frac{\breve{C}_1(a -\phi_{i-k_1,j-k_2;2})(\phi_{i-k_1,j-k_2;1}-\phi_{i,j;1})^2}{\frac{1}{2}Q_2'|\phi_{i,j;2}'|}\\
&+\frac{\breve{C}_2(\phi_{i-k_1,j-k_2;2}-\phi_{i,j;2})(\phi_{i-k_1,j-k_2;1}-\phi_{i,j;1})}{\frac{1}{2}Q_2'|\phi_{i,j;2}'|}\\
&+\frac{\breve{C}_3(\phi_{i-k_1,j-k_2;3}-\phi_{i,j;3})(\phi_{i-k_1,j-k_2;1}-\phi_{i,j;1})}{\frac{1}{2}Q_3'|\phi_{i,j;3}'|}\\
&+\frac{\breve{C}_4(\phi_{i-k_1,j-k_2;2}-\phi_{i,j;2})^2}{\frac{1}{2}Q_2'|\phi_{i,j;2}'|}\\
&+\frac{\breve{C}_5(\phi_{i,j+1;3}-\phi_{i,j;3})(\phi_{i-k_1,j-k_2;2}-\phi_{i,j;2})}{\frac{1}{2}Q_3'|\phi_{i,j;3}'|}\\
&+\frac{\breve{C}_6(a-\phi_{i,j;2})(\phi_{i-k_1,j-k_2;3}-\phi_{i,j;3})^2}{\frac{1}{2}Q_3'|\phi_{i,j;3}'|}\\
\leq&\frac{\breve{C}_1C_2e^{-\eta_2(z_2-k_1 cos\theta -k_2 sin\theta)}[C_2e^{-\eta_2(z_1-k_1cos\theta -k_2 sin\theta)}- C_1
e^{-\eta_2z_1}]^2}{\epsilon_1\rho C_1e^{-\eta_2z_2}}\\
&+\frac{\breve{C}_2[C_2e^{-\eta_2(z_2-k_1cos\theta -k_2 sin\theta)}-C_1e^{-\eta_2z_2}][C_2e^{-\eta_2(z_1-k_1 cos\theta-k_2
sin\theta)}-C_1e^{-\eta_2z_1}]}{\epsilon_1\rho C_1e^{-\eta_2z_2}}\\
&+\frac{\breve{C}_3[C_2e^{\eta_1(z_3-k_1cos\theta -k_2 sin\theta)}-C_1e^{\eta_1z_3}][C_2e^{-\eta_2(z_1-k_1 cos\theta-k_2
sin\theta)}-C_1e^{-\eta_2z_1}]}{\epsilon_1\rho C_1 e^{-\eta_1z_3}}\\
&+\frac{\breve{C}_4[C_2e^{-\eta_2(z_2-k_1cos\theta -k_2 sin\theta)}-C_1 e^{-\eta_2z_2}]^2 }{\epsilon_1\rho C_1 e^{-\eta_2z_2}}\\
&+\frac{\breve{C}_5[C_2e^{\eta_1(z_3-k_1cos\theta -k_2 sin\theta)}-C_1 e^{\eta_1z_3}][C_2e^{-\eta_2(z_2-k_1cos\theta-k_2
sin\theta)}-C_1e^{-\eta_2z_2}]}{\epsilon_1\rho C_1 e^{-\eta_1z_3}}\\
&+\frac{\breve{C}_6C_2 e^{-\eta_2(z_2-k_1 cos\theta -k_2 sin\theta)}[C_2e^{\eta_1(z_3-k_1cos\theta -k_2 sin\theta)}-C_1
e^{\eta_1z_3}]}{\epsilon_1\rho C_1 e^{-\eta_1z_3}}\\
\leq&\tilde{L}_5(e^{-\eta_{2}(z-p_{1}(t))}+e^{-\eta_{2}(z+p_{1}(t))})\leq\tilde{L}_5(e^{\eta_{2}p_{1}(t)}+
e^{\eta_{2}(p_{2}(t)-p_{1}(t))/2})\\
\end{align*}
For $z\geq-p_{2} $, based on Lemma \ref{lem2.1}, Lemma \ref{lem2.2}, Lemma \ref{lem2.3} and \eqref{yin1.4}, we have $z_1=z-p_{1}(t)\geq0$, $z_2=z+p_{1}(t)\geq0$ and
$z_3=z+p_{2}(t)\geq0$. Then it follows  that
\begin{align*}
&|\frac{\mathcal{H}(\phi_{i,j;1},\phi_{i,j;2},\phi_{i,j;3})}{\mathcal{A}(\phi_{i,j;1},\phi_{i,j;2},\phi_{i,j;3})}|\\
\leq &\frac{\dot{C}_1(1-\phi_{i-k_1,j-k_2;3})(\phi_{i-k_1,j-k_2;1}-\phi_{i,j;1})^2}{\frac{1}{2}Q_3'|\phi_{i,j;3}'|}\\
&+\frac{\dot{C}_(1-\phi_{i-k_1,j-k_2;3})(\phi_{i-k_1,j-k_2;2}-\phi_{i,j;2})(\phi_{i-k_1,j-k_2;1}-\phi_{i,j;1})}{\frac{1}{2}Q_3'|\phi_{i,j;3}'|}\\
&+\frac{\dot{C}_3(\phi_{i-k_1,j-k_2;3}-\phi_{i,j;3})(\phi_{i-k_1,j-k_2;1}-\phi_{i,j;1})}{\frac{1}{2}Q_3'|\phi_{i,j;3}'|}\\
&+\frac{\dot{C}_4(1-\phi_{j+1;3})(\phi_{i-k_1,j-k_2;2}-\phi_{i,j;2})^2}{\frac{1}{2}Q_3'|\phi_{i,j;3}'|}\\
&+\frac{\dot{C}_5(\phi_{i-k_1,j-k_2;3}-\phi_{i,j;3})(\phi_{i-k_1,j-k_2;2}-\phi_{i,j;2})}{\frac{1}{2}Q_3'|\phi_{i,j;3}'|}\\
&+\frac{\dot{C}_6(a -\phi_{i,j;2})(\phi_{i-k_1,j-k_2;3}-\phi_{i,j;3})^2}{\frac{1}{2}Q_3'|\phi_{i,j;3}'|}\\
\leq&\frac{\dot{C}_1C_2e^{-\eta_2z_3}[C_2e^{-\eta_2(z_1-k_1cos\theta-k_2 sin\theta)}-C_1 e^{-\eta_2z_1} ]^2}{\epsilon_1\rho
C_1e^{-\eta_2z_3}}\\
&+\frac{\dot{C}_2C_2e^{-\eta_2z_3}[C_2e^{-\eta_2(z_2-k_1 cos\theta-k_2sin\theta)}-C_1e^{-\eta_2z_2}][C_2e^{-\eta_2(z_1-k_1
cos\theta-k_2 sin\theta)}-C_1e^{-\eta_2z_1}]}{\epsilon_1\rho C_1 e^{-\eta_2z_3}}\\
&+\frac{\dot{C}_3[C_2e^{\eta_1(z_3-k_1 cos\theta -k_2 sin\theta)}-C_1 e^{\eta_1z_3}][C_2e^{-\eta_2(z_1-k_1cos\theta-k_2
sin\theta)}-C_1e^{-\eta_2z_1}]}{\epsilon_1\rho C_1 e^{-\eta_2z_3}  }\\
&+\frac{\dot{C}_4C_2e^{-\eta_2z_3}[ C_2e^{-\eta_2(z_2-k_1 cos\theta -k_2 sin\theta)}-C_1 e^{-\eta_2z_2} ]^2}{\epsilon_1\rho C_1
e^{-\eta_2z_3}}\\
&+\frac{\dot{C}_5[C_2e^{-\eta_2(z_3-k_1 cos\theta-k_2 sin\theta)}-C_1 e^{-\eta_2z_3}][C_2e^{-\eta_2(z_2-k_1 cos\theta-k_2
sin\theta)}-C_1e^{-\eta_2z_2}]}{\epsilon_1\rho C_1 e^{-\eta_2z_3}}\\
&+\frac{\dot{C}_6 C_2 e^{-\eta_2z_2}[C_2e^{\eta_1(z_3-k_1 cos\theta -k_2 sin\theta)}-C_1 e^{\eta_1z_3} ]}{\epsilon_1\rho C_1
e^{-\eta_2z_3}}\\
\leq&\tilde{L}_6 (e^{-\eta_{2}(z-p_{1}(t))}+e^{-\eta_{2}(z+p_{1}(t))})\leq\tilde{L}_6(e^{\eta_{2}p_{1}(t)} +
e^{\eta_{2}(p_{2}(t)-p_{1}(t))/2})\\
\end{align*}
Then, we can obtain that
\begin{align*}
|\frac{\mathcal{H}(\phi_{i,j;1}, \phi_{i,j;2},\phi_{i,j;3})}{\mathcal{A}(\phi_{i,j;1}, \phi_{i,j;2},\phi_{i,j;3})}|\leq
L_3(e^{\eta_{2}p_{1}(t)} + e^{\eta_{2}(p_{2}(t)-p_{1}(t))/2}), \; for\; (-p_{2}-p_{1}(t))/2 \leq z\leq-p_{2},
\end{align*}
where $L_3=max\{\tilde{L}_5,\tilde{L}_6\}$.

Let us show that
\begin{align*}
&|\frac{\mathcal{F}(\phi_{i,j;1},\phi_{i,j;2},\phi_{i,j;3})}{ \mathcal{A}(\phi_{i,j;1},\phi_{i,j;2},\phi_{i,j;3})}|\\
&\leq\left\{ {\begin{array}{ll}
L_1'(e^{\eta_{1}p_{1}(t)} + e^{\eta_{1}p_{2}(t)}), \; for\; z\leq0 \mbox \medskip \\
L_2'(e^{\eta_{2}p_{1}(t)} + e^{\eta_{1}(p_{2}(t)-p_{1}(t))/2}), \; for\; 0\leq z\leq(-p_{2}-p_{1}(t))/2 \mbox \medskip \\
L_3'(e^{\eta_{2}p_{1}(t)} + e^{\eta_{2}(p_{2}(t)-p_{1}(t))/2}), \; for\; (-p_{2}-p_{1}(t))/2 \leq z\leq-p_{2}\mbox \medskip \\
\end{array}} \right.
\end{align*}
where $L_1', L_2', L_3'$ are positive constants.

Next, we estimate $|\mathcal{F}/\mathcal{A}|$.
For $z\leq p_1(t)$, we have $z_1=z-p_{1}(t)\leq0$, $z_2=z+p_{1}(t)\leq0$ and $ z_3=z+p_{2}(t)\leq0$, based on Lemma \ref{lem2.1}, Lemma \ref{lem2.2}{\rm(ii)}, Lemma \ref{lem2.3} and \eqref{yin1.1}. Then it follows  that
\begin{align*}
&|\frac{\mathcal{F}(\phi_{i,j;1},\phi_{i,j;2},\phi_{i,j;3})}{ \mathcal{A}(\phi_{i,j;1},\phi_{i,j;2},\phi_{i,j;3})}|\\
\leq&|\frac{2(1-\phi_{i,j;1})\phi_{i,j;2}[\phi_{i,j;1}H_{11}+(\phi_{i,j;3}-a)H_{11}]}{Q_1'|\phi_{i,j;1}'|}|\leq \frac{2
\frac{1}{\rho}|\phi_{i,j;2}|[C+ \frac{1}{\rho}|\phi_{i,j;3}'|C]}{Q_1'}\\
\leq&\frac{2 C|\phi_{i,j;2}'|+2C|\phi_{i,j;3}'|}{\epsilon_1\rho^2}\leq\frac{2 C C_0((e^{\eta_{1}(z+p_{1}(t))} +
e^{\eta_{1}(z+p_{2}(t))}))}{\epsilon_1\rho^2}\\
\leq&\bar{L}_1(e^{\eta_{1}(z+p_{1}(t))}+e^{\eta_{1}(z+p_{2}(t))})\leq\bar{L}_1(e^{\eta_{1}p_{1}(t)}+e^{\eta_{1}p_{2}(t)})\\
\end{align*}

For $p_1(t)\leq z\leq0$,we have $z_1=z-p_{1}(t)\geq0$, $z_2=z+p_{1}(t)\leq0$ and $ z_3=z+p_{2}(t)\leq0$, based on Lemma \ref{lem2.1}, Lemma \ref{lem2.2}{\rm(ii)}, Lemma \ref{lem2.3},
and \eqref{yin1.2}. Then it
follows that
\begin{align*}
&|\frac{\mathcal{F}(\phi_{i,j;1},\phi_{i,j;2},\phi_{i,j;3})}{ \mathcal{A}(\phi_{i,j;1},\phi_{i,j;2},\phi_{i,j;3})}|\\
\leq&|\frac{2(1-\phi_{i,j;1})\phi_{i,j;2}[\phi_{i,j;1}H_{11}+(\phi_{i,j;3}-a)H_{11}]}{Q_1'|\phi_{i,j;1}'|+Q_2'|\phi_{i,j;2}'|}|\\
\leq&\frac{2|(1-\phi_{i,j;1})||\phi_{i,j;2}||\phi_{i,j;1}||H_{11}|}{Q_1'|\phi_{i,j;1}'|}+
\frac{2|(1-\phi_{i,j;1})||\phi_{i,j;2}||\phi_{j;3}-a)||H_{11}|}{Q_2'|\phi_{i,j;1}'|}\\
\leq&\frac{2C|\phi_{i,j;2}'|}{\epsilon_1\rho^2}+\frac{2C|\phi_{i,j;3}'|}{\epsilon_2\rho^2} \leq
\frac{2C}{\epsilon_1\rho^2}e^{\eta_{1}(z+p_{1}(t))}+\frac{2C}{\epsilon_2\rho^2}e^{\eta_{1}(z+p_{2}(t))}\\
\leq&\bar{L}_2(e^{\eta_{1}(z+p_{1}(t))}+e^{\eta_{1}(z+p_{2}(t))})\leq\bar{L}_2 (e^{\eta_{1}p_{1}(t)} + e^{\eta_{1}p_{2}(t)})\\
\end{align*}

We can obtain that
\begin{align*}
|\frac{\mathcal{F}(\phi_{i,j;1},\phi_{i,j;2},\phi_{i,j;3})}{\mathcal{A}(\phi_{i,j;1},\phi_{i,j;2},\phi_{i,j;3})}|\leq
L_1'(e^{\eta_{1}p_{1}(t)} + e^{\eta_{1}p_{2}(t)}), \; for\; z\leq0,
\end{align*}
where $L_1'=max\{\bar{L}_1,\bar{L}_2\}$.

For $0\leq z\leq-p_{1}(t)$, we have $z_1=z-p_{1}(t)\geq0$, $z_2=z+p_{1}(t)\leq0$ and $z_3=z+p_{2}(t)\leq0$, based on Lemma \ref{lem2.1}, Lemma \ref{lem2.3}, \eqref{2.26} and \eqref{yin1.2}. Then it follows that
\begin{align*}
&|\frac{\mathcal{F}(\phi_{i,j;1},\phi_{i,j;2},\phi_{i,j;3})}{ \mathcal{A}(\phi_{i,j;1},\phi_{i,j;2},\phi_{i,j;3})}|\\
&\leq|\frac{2(1-\phi_{i,j;1})\phi_{i,j;2}[\phi_{i,j;1}H_{11}+(\phi_{j;3}-a)H_{12}]}{Q_2'|\phi_{i,j;2}'|}|\\
&\leq\frac{2|(1-\phi_{i,j;1})||\phi_{i,j;2}||\phi_{i,j;1}||H_{11}|}{Q_2'|\phi_{i,j,2}'|}+\frac{2|(1-\phi_{i,j;1})||\phi_{i,j;2}||(\phi_{i,j;3}-a)||H_{12}|}{Q_2'|\phi_{i,j;2}'|}\\
&\leq\bar{L}_3 (e^{-\eta_{2}(z-p_{1}(t))}+e^{\eta_{1}(z+p_{2}(t))})\leq\bar{L}_3 (e^{\eta_{2}p_{1}(t)}+
e^{\eta_{1}(p_{2}(t)-p_{1}(t))/2})\\
\end{align*}

For $-p_{1}(t)\leq z\leq(-p_{2}-p_{1}(t))/2 $,
we have $z_1=z-p_{1}(t)\geq0$, $z_2=z+p_{1}(t)\geq0$ and $ z_3=z+p_{2}(t)\leq0$, based on Lemma \ref{lem2.1}, Lemma \ref{lem2.2}{\rm(ii)}, Lemma \ref{lem2.3}, and \eqref{yin1.3}.
Then it follows that
\begin{align*}
&|\frac{\mathcal{F}(\phi_{i,j;1},\phi_{i,j;2},\phi_{i,j;3})}{ \mathcal{A}(\phi_{i,j;1},\phi_{i,j;2},\phi_{i,j;3})}|\\
\leq&|\frac{2(1-\phi_{i,j;3})(a-\phi_{i,j;2})[\phi_{i,j;1}H_{21}+(\phi_{i,j;3}-a)H_{22}]}{Q_2'|\phi_{i,j;2}'|}|\\
\leq&\frac{2|(1-\phi_{i,j;3})||a-\phi_{i,j;2}||\phi_{i,j;1}||H_{11}|}{Q_2'|\phi_{i,j;2}'|}+\frac{2|(1-\phi_{i,j;3})||a-\phi_{i,j;2}||(\phi_{i,j;3}-a)||H_{22}|}{Q_2'|\phi_{i,j;2}'|}\\
\leq&\bar{L}_4(e^{-\eta_{2}(z-p_{1}(t))}+e^{\eta_{1}(z+p_{2}(t))})\leq\bar{L}_4 (e^{\eta_{2}p_{1}(t)} +
e^{\eta_{1}(p_{2}(t)-p_{1}(t))/2})\\
\end{align*}

We can obtain that
 \begin{align*}
|\frac{\mathcal{F}(\phi_{i,j;1},\phi_{i,j;2},\phi_{i,j;3})}{ \mathcal{A}(\phi_{i,j;1},\phi_{i,j;2},\phi_{i,j;3})}|\leq
L_2'(e^{\eta_{2}p_{1}(t)} + e^{\eta_{1}(p_{2}(t)-p_{1}(t))/2}), \; for\; 0\leq z\leq(-p_{2}-p_{1}(t))/2,
\end{align*}
where $L_2'=max\{\bar{L}_3,\bar{L}_4\}$.

For $(-p_{2}-p_{1}(t))/2\leq z\leq-p_{2}$, we have $z_1=z-p_{1}(t)\geq0$, $z_2=z+p_{1}(t)\geq0$ and $z_3=z+p_{2}(t)\leq0$,
based on Lemma \ref{lem2.1}, Lemma \ref{lem2.2}{\rm(ii)}, Lemma \ref{lem2.3}, and \eqref{yin1.3}.
Then it follows that
\begin{align*}
&|\frac{\mathcal{F}(\phi_{i,j;1},\phi_{i,j;2},\phi_{i,j;3})}{ \mathcal{A}(\phi_{i,j;1},\phi_{i,j;2},\phi_{i,j;3})}|\\
\leq&|\frac{2(1-\phi_{i,j;3})(a-\phi_{i,j;2})[\phi_{i,j;1}H_{21}+(\phi_{i,j;3}-a)H_{22}]}{Q_2'|\phi_{i,j;2}'|+Q_3'|\phi_{i,j;3}'|}\\
\leq&\frac{2|(1-\phi_{i,j;3})||a-\phi_{i,j;2}||\phi_{i,j;1}||H_{11}|}{Q_2'|\phi_{i,j;2}'|}+\frac{2|(1-\phi_{i,j;3})||a-\phi_{i,j;2}|
|(\phi_{i,j;3}-a)||H_{22}|}{Q_3'|\phi_{i,j;3}'|}\\
\leq& \bar{L}_5(e^{-\eta_{2}(z-p_{1}(t))}+e^{-\eta_{2}(z+p_{1}(t))})\leq\bar{L}_5(e^{\eta_{2}p_{1}(t)}+
e^{\eta_{2}(p_{2}(t)-p_{1}(t))/2})\\
\end{align*}

For $z\geq-p_{2}$, we have $z_1=z-p_{1}(t)\geq0$, $z_2=z+p_{1}(t)\geq0$ and $z_3=z+p_{2}(t)\geq0$,
based on Lemma \ref{lem2.1}, Lemma \ref{lem2.2}{\rm(ii)}, Lemma \ref{lem2.3}, and \eqref{yin1.4}.
Then it follows that
\begin{align*}
&|\frac{\mathcal{F}(\phi_{i,j;1},\phi_{i,j;2},\phi_{i,j;3})}{ \mathcal{A}(\phi_{i,j;1},\phi_{i,j;2},\phi_{i,j;3})}|\\
\leq&|\frac{2(1-\phi_{i,j;3})(a-\phi_{i,j;2})[\phi_{i,j;1}H_{21}+(\phi_{i,j;3}-a)H_{22}]}{Q_3'|\phi_{i,j;3}'|}\\
\leq&\frac{2}{\epsilon_3\rho}[C|\phi_{i,j;1}|+C|a-\phi_{i,j;2}|]\leq\frac{2C}{\epsilon_3\rho^2}(e^{\eta_{1}(z-p_{1}(t))}+e^{\eta_{1}(z+p_{1}(t))})\\
\leq&\bar{L}_6 (e^{-\eta_{2}(z-p_{1}(t))} + e^{-\eta_{2}(z+p_{1}(t))})\leq\bar{L}_6 (e^{\eta_{2}p_{1}(t)} +
e^{\eta_{2}(p_{2}(t)-p_{1}(t))/2})\\
\end{align*}

We can obtain that
 \begin{align*}
|\frac{\mathcal{F}(\phi_{i,j;1},\phi_{i,j;2},\phi_{i,j;3})}{ \mathcal{A}(\phi_{i,j;1},\phi_{i,j;2},\phi_{i,j;3})}|\leq
L_3'(e^{\eta_{2}p_{1}(t)} + e^{\eta_{2}(p_{2}(t)-p_{1}(t))/2}), \; for\; (-p_{2}-p_{1}(t))/2\leq z\leq-p_{2},
\end{align*}
where $L_3'=max\{\bar{L}_5,\bar{L}_6\}$. For all $t\leq t_0$, where $M_1=max\{L_1,L_2,L_3\}$  and  $M_2=max\{L_1',L_2',L_3'\}$.
Then choosing $L>(M_1,M_2)$, it follows that
\begin{align*}
&F(\overline{U})(i,j,t)=-Le^{\kappa p_{1}}Q_1'\phi_{i,j;1}' +Le^{\kappa p_{1}}Q_2' \phi_{i,j;2}'+Le^{\kappa p_{1}}Q_3' \phi_{i,j;3}'-
\mathcal{F}(\phi_{i,j;1},\phi_{i,j;2},\phi_{i,j;3})\\
&-\mathcal{H}(\phi_{i,j;1},\phi_{i,j;2},\phi_{i,j;3})
\geq \mathcal{A}(\phi_{i,j;1},\phi_{i,j;2},\phi_{i,j;3})[L-(M_1+M_2)]e^{\kappa p_1(t)}\geq 0.\\
\end{align*}
By \eqref{eq2.11}, \eqref{eq2.12}, Lemma \ref{lem2.4}, hence $\overline{U}_{i,j}(t)$ is a super-solution of \eqref{eq1.1} for $t<t_0$.
Similarly, we can prove that $\underline{U}_{i,j}(t)$ is a subsolution of \eqref{eq1.1} for $t<t_0$.
Then we obtain
\begin{align*}
&\overline{U}_{i,j}(t)-\underline{U}_{i,j}(t)\\
=&Q\big(\phi_{i,j;1}(j+\bar{c}t-p_{1}(t)),\phi_{i,j;2}(j+\bar{c}t+p_{1}(t)),\phi_{i,j;3}(j+\bar{c}t+p_{2}(t))\big)\\
&-Q\big(\phi_{i,j;1}(j+\bar{c}t-r_{1}(t)),\phi_{i,j;2}(j+\bar{c}t+r_{1}(t)),\phi_{i,j;3}(j+\bar{c}t+r_{2}(t))\big)\\
=&\int_0^1F(\phi_{i,j;1}(j-\theta p_1-(1-\theta) r_{1}),\phi_{i,j;2}(j+\theta p_1+(1-\theta) r_{1}),\phi_{i,j;3}(j+\theta
p_2+(1-\theta) r_{2}))d\theta\\
&\times(p_1-r_1).
\end{align*}
We have
\begin{align*}
\overline{U}_{i,j}(t) >\underline{U}_{i,j}(t)\ \ \mbox{and}\ \
\sup_{i,j\in\mathbb{Z}}(\overline{U}_{i,j}(t)-\underline{U}_{i,j}(t))
\leq Ce^{ \kappa s_{1}t},~~\mbox{for }t\leq0
\end{align*}
hence the lemma is proved.

From \cite{Dongli}, existence and uniqueness of entire solution of \eqref {eq1.1} can be shown.
There exists a unique entire solution $u_{i,j}(t)$ of \eqref {eq1.1} such that
$ \underline{U}_{i,j}(t) \leq u_{i,j}(t)\leq  \overline{U}_{i,j}(t)$
for $t<T$.
Then we Define
\begin{align}
\omega:=-\frac{1}{\kappa}log(e^{-\kappa r_0}-\frac{L}{s_1}).
\end{align}
Furthermore, on the premise that the supersolution and subsolution are established. By \cite[Theorem 3.3]{ZhangLiWu}, we can get the asymptotic behavior of
\eqref{eq1.1}.
\section{Proof of theorem 1.3}

Since the proof of theorem 1.3 is quite similar to that of theorem 1.2, we only point out the main differences in this section. Similar to that of theorem 1.2.
Set $c_1:=\hat{v}$, $c_2:=\tilde{v_1}$ and $c_3:=\hat{v_1}$. Let
$\phi_{i,j;1}=\hat\varphi_{i,j;1}$, $\phi_{i,j;2}=\tilde{\varphi}_{i,j;2}$, $\phi_{i,j;3}=\hat{\varphi}_{i,j;2}$, $j\in \mathbb{Z}, k=1,2,3,$
 be traveling fronts of \eqref{eq1.1} that satisfy

\begin{align*}
\left\{
\begin{array}{ll}
{c_k}\phi_{i,j;k}'(\xi)=\sum_{k_1}\sum_{k_2}J(k_1,k_2)\phi_{i-k_1,j-k_2;k}(\xi)-\phi_{i,j;k}(\xi)+f_{i,j}(\phi_{i,j;k}(\xi)),\medskip\\
\phi_{i+N_1,j;k}(\xi)=\phi_{i,j;k}(\xi)=\phi_{i,j+N_2;k}(\xi),    i,j\in \mathbb{Z},~\xi\in \mathbb{R}.\\
 \phi_{i,j;k}(-\infty)=\alpha_i,\phi_{i,j;k}(\infty)=\beta_i
\end{array}\right.
\end{align*}
where $(\alpha_1,\beta_1,\alpha_2,\beta_2,\alpha_3,\beta_3)=(1,0,0,a,a,0)$.

First, we consider the auxiliary rational function $\tilde Q(y,z,w)$
\begin{align*}
\tilde Q(y,z,w):=&z+\frac{(1- y )z(a-w)(-z)+y(a-z)w(1-z)}{(1-y)za+(a-z)w}.
\end{align*}
Then we take the function $\tilde p_{i}'(t)$ and $\tilde r_{i}'(t)$, $i=1,2$ are the solutions of the following initial value problems:
\begin{align}
\tilde p_{1}'(t)=s_1+Le^{\kappa \tilde p_{1}(t)},\,\,\tilde p_{1}(0)=\tilde p_0;
\label{eq3.1}\\
\tilde r_{1}'(t)=s_1-Le^{\kappa \tilde r_{1}(t)},\,\,\tilde r_{1}(0)=\tilde r_0;
\label{eq3.2}\\
\tilde p_{2}'(t)=s_2-Le^{\kappa \tilde p_{1}(t)},\,\,\tilde p_{2}(0)=\tilde r_0;
\label{eq3.3}\\
\tilde r_{2}'(t)=s_2+Le^{\kappa \tilde r_{1}(t)},\,\,\tilde r_{2}(0)=\tilde p_0,
\label{eq3.4}
\end{align}
where $L,s_1,s_2,\tilde p_0,\tilde r_0$ and $\kappa$ are the same as in Section 2.

Define the function $\overline{U}^{*}_{i,j}(t)$ and $\underline{U}^{*}_{i,j}(t)$ as follows.
\begin{align*}
&\overline{U}^{*}_{i,j}(t)=\tilde Q\big(\phi_{i,j;1}(j+\bar{c}t-\tilde p_{1}(t)),\phi_{i,j;2}(j+\bar{c}t+\tilde
p_{1}(t)),\phi_{i,j;3}(j+\bar{c}t+\tilde p_{2}(t))\big),\\
&\underline{U}^{*}_{i,j}(t)=\tilde Q\big(\phi_{i,j;1}(j+\bar{c}t-\tilde r_{1}(t)),\phi_{i,j;2}(j+\bar{c}t+\tilde
r_{1}(t)),\phi_{i,j;3}(j+\bar{c}t+\tilde r_{2}(t))\big)            
\end{align*}

By the similar argument as in the function $\overline{U}^{*}_{i,j}(t)$ and $\underline{U}^{*}_{i,j}(t)$ are a pair of super-sub-solution of
\eqref{eq1.1} for $t\leq t_0$ with some constant $t_0<0$.
Hence, the existence and uniqueness of the entire solution $u_{i,j}(t)$ of \eqref{eq1.1} can be shown and $u_{i,j}(t)$ satisfies
\begin{align*}
\underline{U}^{*}_{i,j}(t)\leq u_{i,j}(t)\leq\overline{U}^{*}_{i,j}(t)
\end{align*}
for all $i,j\in Z$ and $t\leq t_0$, Moreover, it is not difficult to check the entire solution  $u_{i,j}(t)$ satisfies by using the
similar argument as in and taking
\begin{align*}
\omega_1:=-\frac{1}{\kappa}log(e^{-\kappa \tilde r_{0}}-\frac{L}{s_1}),\omega_2:=-\frac{1}{\kappa}log(e^{-\kappa \tilde
r_{0}}-\frac{L}{s_1})+\tilde p_{0}+\tilde r_{0},
\end{align*}
Similar to the proof argument as Section 2, we can also obtain that $\overline{U}^{*}_{i,j}(t)$ and $\underline{U}^{*}_{i,j}(t)$
are a pair of supersolution and subsolution of \eqref{eq1.1}.
Since the proof of Theorem \ref{thm1.3} is quite similar to that of theorem \ref{thm1.2}, we omit the detail of the proof.


\begin{thebibliography}{99}
\bibitem{ChenGuoNinomiyaYao}Y.-Y. Chen, J.-S. Guo, H. Ninomiya and C.-H. Yao, Entire solutions with merging three fronts to the
    Allen-Cahn equation, mathscidoc:1609.03007.
\bibitem{ZhangLiWu} S. L. Wu and C. H. Hsu, Entire solutions with merging fronts to a bistable periodic lattice dynamical system,
    Discrete Contin. Dyn. Syst., 36 (2016), 2329--2346.
\bibitem{YYC} Y.-Y. Chen, Entire solution originating from three fronts for a discrete diffusive equation, Tamkang Journal of
    Mathematics 48 (2017), 215--226.
\bibitem{CWG} C. P. Cheng, W. T. Li and G. Lin, Travelling wave solutions in periodic monostable equations on a two-dimensional spatial lattice,
 IMA J. Appl. Math., 80 (2015), 1254--1272.
\bibitem{JSGCHWE}J. S. Guo and C. H. Wu, Existence and uniqueness of traveling waves for a monostable 2-D lattice dynamical system, Osaka J. Math., 45 (2008), 327--346.

\bibitem{JSGCHWT} J. S. Guo and C. H. Wu, Traveling wave front for a two-component lattice dynamical system arising in competition models, J. Differential Equations, 252 (2012), 4357--4391.

\bibitem{SX}S. Ma and X. Zou, Propagation and its failure in a lattice delayed differential equation with global interaction, J.
    Differential Equations, 212 (2005), 129--190.
\bibitem{CCWU}C. C. Wu, Uniqueness of traveling waves for a two-dimensional bistable periodic lattice dynamical system, Abstr. Appl.
    Anal., Volume 2012, Article ID 289168, 10 pages.
\bibitem{Dongli} F.D. Dong, W.T. Li and L.Zhang,Entire solutions in a two-dimensional nonlocal lattice dynamical system,{\it Comm.
    Pure Appl. Anal.}, 17(2018), 2517--2545.

\bibitem{JCWU}J. S. Guo and C. H. Wu, Entire solutions for a two-component competition system in a lattice, Tohoku Math. J., 62 (2010), 17--28.

\bibitem{CHWU} C. H. Wu, A general approach to the asymptotic behavior of traveling waves in a class of three-component lattice dynamical systems, J. Dynam. Differential Equations, 28 (2016), 317--338.
\bibitem{Wang7} Z.C. Wang, W.T. Li and S. Ruan, Entire solutions in delayed lattice differential equations with
monostable nonlinearity, {\it SIAM J. Math. Anal.}, 40 (2009), 2392--2420.
\bibitem{wushiyang}	S.L. Wu, Z.X. Shi and F.Y. Yang, Entire solutions in periodic lattice dynamical systems, {\it J. Differential
    Equations},  255  (2013), 3505--3535.
\bibitem{Moritan} Y. Morita and H. Ninomiya, Entire solutions with merging fronts to reaction-diffusion equations,
{\it J. Dynam. Differential Equations}, 18 (2006), 841--861.
\bibitem{Guom} J.-S. Guo and Y. Morita, Entire solutions of reaction-diffusion equations and an application
to discrete diffusive equations, {\it Discrete Contin. Dyn. Syst.}, 12 (2005), 193--212.
\bibitem{ChenGuo}X. Chen, J. S. Guo and C. C. Wu, Traveling waves in discrete periodic media for bistable dynamics, Arch. Ration.
    Mech. Anal., 189 (2008), 189--236.
\bibitem{JSGFH} J. S. Guo and F. Hamel, Front propagation for discrete periodic monostable equations, Math. Ann., 335 (2006), 489--525.
\bibitem{JSGCHW}  J. S. Guo and C. H. Wu, Front propagation for a two-dimensional periodic monostable lattice dynamical system, Discrete Contin. Dyn.Syst., 26 (2010), 197--223.
\bibitem{SMPW} S. Ma, P. Weng and X. Zou, Asymptotic speed of propagation and traveling wavefronts in a non-local delayed lattice differential equation, Nonlinear Anal., 65 (2006), 1858--1890.
\bibitem{PWAC}P. W. Bates and A. Chmaj, A discrete convolution model for phase transitions, Arch. Ration. Mech. Anal., 150 (1999), 281--305.
\bibitem{CPWT}C. P. Cheng, W. T. Li and G. Lin, Travelling wave solutions in periodic monostable equations on a two-dimensional spatial lattice, IMA J. Appl. Math., 80 (2015), 1254--1272.
 \bibitem{CPYH}C. P. Cheng, Y. H. Su and Z. Feng, Wave propagation for monostable 2-D lattice differential equations with delay, Internat. J.Bifur. Chaos Appl. Sci. Engrg., 23 (2013), 1350077, 11 pp.
\bibitem{CPWTZC} C. P. Cheng, W. T. Li and Z. C. Wang, Persistence of bistable waves in a delayed population model with stage structure on a two-dimensional spatial lattice, Nonlinear Anal. Real World Appl., 13 (2012), 1873--1890.
\bibitem{CPCWT} C. P. Cheng, W. T. Li and Z. C. Wang, Asymptotic stability of traveling wavefronts in a delayed population model with stage structure on a two-dimensional spatial lattice, Discrete Contin. Dyn. Syst. Ser. B, 13 (2010), 559--575.
\bibitem{CPCWTZ} C. P. Cheng, W. T. Li and Z. C. Wang, Spreading speeds and travelling waves in a delayed population model with stage structure on a 2D spatial lattice, IMA J. Appl. Math., 73 (2008), 592--618.
\bibitem{JACU} J. Carr and A. Chmaj, Uniqueness of travelling waves for nonlocal monostable equations, Proc. Amer. Math. Soc., 132 (2004), 2433--2439.
\bibitem{XCSC} X. Chen, S. C. Fu and J. S. Guo, Uniqueness and asymptotics of traveling waves of monostable dynamics on lattices, SIAM J. Math. Anal., 38 (2006), 233--258.
\bibitem{XCJSU} X. Chen and J. S. Guo, Uniqueness and existence of traveling waves for discrete quasilinear monostable dynamics, Math. Ann., 326 (2003), 123--146.

\bibitem{JSGYW} J. S. Guo, Y. Wang, C. H. Wu and C. C. Wu, The minimal speed of traveling wave solutions for a diffusive three species competition system, Taiwanese J. Math., 19 (2015), 1805--1829.

\end{thebibliography}
\end{document}